\documentclass[onefignum,onetabnum]{siamart190516}

\usepackage{listofitems}
\usepackage{graphicx}
\usepackage{comment}
\excludecomment{OldVersion}
\usepackage{bbm}
\usepackage{stmaryrd}
\usepackage{pgfplotstable}
\usepackage{pgfplots}
\usepackage{tikz}
\usepackage{tikz-qtree}
\usepackage{algorithm}
\usepackage{algpseudocode}
\usepackage{amsmath,amsfonts}
\usepackage{caption}
\usepackage{subcaption}

\usepackage{amssymb}
\usepackage[default]{mycolor}

\newenvironment{compositeproof}{\begin{proof}}{\end{proof}}

\let\Halmos\

\newcommand{\QED}{\hfill \smash{\scalebox{.65}[1]{\ensuremath{\Box}}}}

\headers{MIP relaxations of composite functions}{Taotao He, Mohit Tawarmalani}

\newsiamremark{remark}{Remark}
\newsiamremark{hypothesis}{Hypothesis}
\crefname{hypothesis}{Hypothesis}{Hypotheses}
\newsiamremark{example}{Example}

\DeclareMathOperator{\proj}{proj}

\DeclareMathOperator{\RLT}{RLT}
\DeclareMathOperator{\conv}{conv}
\DeclareMathOperator{\conc}{conc}

\DeclareMathOperator{\graph}{gr}

\DeclareMathOperator{\aff}{aff}

\DeclareMathOperator{\vertex}{vert}

\DeclareMathOperator{\hypo}{hyp}

\DeclareMathOperator{\R}{\mathbb{R}}
\DeclareMathOperator{\Z}{\mathbb{Z}}

\DeclareMathOperator{\ri}{ri}

\newcommand{\revision}{}

\title{MIP relaxations in factorable programming\thanks{Submitted to the editors DATE.
\funding{The work was supported in part by NSFC Grant 21Y1254213, the Shanghai Pujiang Grant 21PJC073, and by NSF CMMI Grant 1727989}}}

\author{Taotao He\thanks{Antai College of Management and Economics, Shanghai Jiao Tong University, \email{hetaotao@sjtu.edu.cn}}
\and Mohit Tawarmalani\thanks{Krannert School of Management, Purdue University, \email{mtawarma@purdue.edu}.}}

\begin{document}

%------- shortcuts-------------
\def \mcirc {\mathop{\circ}}

\def \CR {R} % composite relaxation
\def \J {{\bar{J}}} %complement of subset J
\def \s {{\bar{s}}} %fixed s
\def \u {{\bar{u}}} %fixed u
\def \p {{\bar{\phi}}} %fixed phi
\def \z {{\bar{z}}} %fixed z
\def \ts {{\tilde{s}}}
\def \tu {{\tilde{u}}}
\def \ta {{\tilde{\alpha}}}
\def \tb {{\tilde{b}}}
\def \mJ {{\mathcal{J}}} %collection of indexes
\def \mX {{\mathcal{X}}}
\def \dbl {{[\![}}
\def \dbr {{]\!]}}
%collection of indexes

\maketitle

% REQUIRED
\begin{abstract}
In this paper, we develop new discrete relaxations for nonlinear expressions in factorable programming. We utilize specialized convexification results as well as composite relaxations to develop mixed-integer programming \revision{(MIP)} relaxations. Our relaxations rely on ideal formulations of convex hulls of outer-functions over a combinatorial structure that captures local inner-function structure. The resulting relaxations often require fewer variables and are tighter than currently prevalent ones. Finally, we provide computational evidence to demonstrate that our relaxations close approximately 60-70\% of the gap relative to McCormick relaxations and significantly improves the relaxations used in a state-of-the-art solver on various instances involving polynomial functions.
\end{abstract}

% REQUIRED
\begin{keywords}
  Mixed-Integer Nonlinear Programming, MIP Relaxation, Convexification, Incremental Formulation, Relaxation Propagation
\end{keywords}

% REQUIRED
\begin{AMS}
  68Q25, 68R10, 68U05
\end{AMS}

\def \ephi {\bar{\phi} }
\section{Introduction}

Mixed-integer nonlinear programming (MINLP) algorithms and software rely primarily on three techniques for constructing relaxations. First, factorable programming techniques are used to relax composite functions 
\cite{mccormick1976computability,tawarmalani2004global,misener2014antigone,mahajan2021minotaur,bestuzheva2023global}. Second, there has been a growing interest in MIP relaxations motivated by significant advances in solver capability for this class of problems~\cite{geissler2011using,misener2012global,geissler2011using,gupte2013solving,huchette2019combinatorial,nagarajan2019adaptive,beach2022compact}. Third, solvers often utilize envelope results for special functions over hypercubes to tighten relaxations; see, for example, \cite{rikun1997convex,sherali1997convex,tawarmalani2002convexification,cafieri2010convex,anstreicher2010computable,tawarmalani2013explicit,khajavirad2013convex,gupte2020extended}. Regardless, solvers arguably ignore much of inner-function structure while relaxing composite functions, particularly while constructing discretized relaxations or if a specialized result is not already available~\cite{tawarmalani2004global,misener2014antigone,nagarajan2019adaptive,mahajan2021minotaur,bestuzheva2023global}. Recently, \cite{he2021new,he2022tractable} has devised techniques to exploit the inner-function structure for constructing convex relaxations. Even so, the typical strategy of constructing MIP relaxations is to decompose composite functions into univariate and multilinear functions, and then to focus on discretizing the multilinear expressions. As yet, the literature does not provide sufficient guidance on how to use (i) specialized convexification results and/or (ii) composite relaxations so that they are not agnostic of the discretization strategy.

%Classical factorable programming does not utilize any structural information on the inner function besides bounds~\cite{mccormick1976computability}, while special-structured results are only available for specific function types~\cite{rikun1997convex,sherali1997convex,benson2004concave,meyer2004trilinear}. The approach introduced recently in \cite{he2021new,he2022tractable} develops a framework that lies between these extremes. The technique allows the inner-function structure to be exploited in a completely generic fashion without the need for special-structure identification \cite{he2021new} and develops tractable relaxations when outer-function satisfies certain technical conditions \cite{he2022tractable}.

In this paper, we show that there are generic ways to construct MIP relaxations that address both these challenges.
%[82, 69, 52, 171] that incorporate many advanced techniques [1, 2, 25, 23, 92, 112, 24]
First, we show that envelope results naturally yield ideal formulations for discretized functions.  Second, we show that our constructions seamlessly integrate techniques introduced in \cite{he2021new}, thus leveraging inner-function structure in MIP relaxations for composite functions. 
%Towards this end, we discover an affine transformation that relates the simplotope used to discretize a function and the one encodes the inner-function structure.
\revision{Towards this end, we discover an affine map between the simplotope used to discretize a function and the one used to encode the inner-function structure.}
%interpret a specific simplotope as a projection of another. The projected simplotope is used to discretize the function, while the extended simplotope is derived from a polytope that encodes the inner-function structure. 
The former allows for use of envelope results to produce ideal formulations while the latter captures the inner-function structure. 
%The latter simplotope allows for use of convex/concave envelope results to capture the inner-function structure while the former is used to produce ideal relaxations. 
%The envelope results that are useful in this context are those obtained over a simplotope instead of a hypercube.

Addressing these challenges leads to several improvements. The resulting relaxations are tighter and often more economical than prevalent MIP relaxations. For example, when the envelopes are polynomially separable, our relaxations require exponentially fewer continuous variables. Interestingly, ideality is guaranteed by a generic technique we use to derive envelope results. Finally, our relaxations can, upon branching, exploit local inner-function structure to further reduce the relaxation gap. \revision{The current state-of-the-art solvers \cite{tawarmalani2004global,mahajan2021minotaur,bestuzheva2023global} use spatial branching to improve relaxation quality. We show that, on some polynomial optimization problems, our MIP (resp. LP) relaxations take almost the same time (resp. significantly less time) to solve but close up to 50\% (resp. up to 40\%) of the gap relative to relaxation in \texttt{SCIP} derived at the root node after range-reduction and cutting planes. Moreover, in Section~\ref{section:comparison}, we show how our MIP relaxation can be tightened using local bounds on estimators obtained via spatial branching. }
%so that, upon branching, they get tighter by automatically update local bounds for underestimators of inner-functions. 

\revision{
The paper is laid out as follows. In Section~\ref{section:MIP-conv}, we construct ideal MIP relaxations for discretized outer-functions. To do so, we utilize convex/concave envelopes alongside an incremental formulation that models the union of certain faces of a specific simplotope. In Section~\ref{section:MIP-relaxation}, we integrate these MIP formulations with composite relaxations~\cite{he2021new}. In Section~\ref{section:comparison}, we develop geometric insights into the quality of relaxations and discuss ways to tighten the relaxations. Finally, in Section~\ref{section:computation}, we develop a propagation strategy, Algorithm~\ref{alg:propagation}, for developing MIP relaxations of expression trees, and show that our relaxations close significant gap on problems relative to state-of-the-art solvers without incurring significant computational overhead. Last, in Section~\ref{section:extension}, we discuss logarithmic MIP relaxations, and MIP relaxations for multiple composite functions. 
}

\paragraph{\bf{Notation}} 
%Throughout this paper we use the following notation.
We shall denote the convex hull of set $S$ by $\conv(S)$, the projection of a set $S$ to the space of $x$ variables by $\proj_x(S)$, the extreme points of $S$ by $\vertex(S)$, and the convex (resp. concave) envelope of $f(\cdot)$ over $S$ by $\conv_{S}(f)(\cdot)$ (resp. $\conc_S(f)(\cdot)$). For a vector of functions $g:D \to \R^m$, its hypograph of $g$ is defined as $\hypo(g): =\{ (x,\mu) \mid \mu \leq g(x), x \in D \}$, and its graph is defined as $\graph(g): =\{ (x,\mu) \mid \mu =  g(x), x \in D \}$. We  denote an index set $\{1, \ldots,n\}$ by $[ n ]$, where $n$ is a positive integer.
%\section{Preliminaries}\label{section:exp}
%\def\L{\mathbb{L}}
%\def\G{\mathcal{G}}
%\def\extreme{\mathop{\text{ext}}}
%\def\Qsimplex{\upsilon}
\def\V{\mathcal{V}}
\def\sfunct{\eta}
\def\DO{\mathcal{D}}
\def\BO{\mathcal{B}}
\def\mcirc{\mathop{\circ}}

\section{Preliminaries}\label{section:discrete-relaxations}
Consider a composite function $\phi \mcirc f \colon X \subseteq \R^m \to \R$ defined as $(\phi \mcirc f)(x) = \phi\bigl(f(x)\bigr)$, where $f(x) := \bigl(f_1(x), \ldots, f_d(x) \bigr)$. We refer to each $f_i \colon \R^m \to \R$  as an \textit{inner-function} while $\phi \colon \R^d \to \R$ is referred to as the \textit{outer-function}. We assume that the inner-function is bounded, that is, there exist two vectors $f^L$ and $f^U$ such that $f^L\leq f(x) \leq f^U$ for every $x \in X$. In this section, we present definitions relevant to  MIP relaxations, and a discretization scheme for the hypograph of $\phi \mcirc f$. 

 For a set $x \in S \subseteq \R^{n_1}$, let $R \subset \R^{n_1}\times\R^{n_2}\times\R^{n_3}$ and define $E:=\bigl\{(x,y,\delta)\in \R^{n_1}\times\R^{n_2}\times\{0,1\}^{n_3} \bigm| (x,y,\delta) \in R \bigr\}$. If $S\subseteq \proj_{x}(E)$, we refer to $E$ as an \textit{MIP relaxation} of $S$  and $R$ is the \textit{continuous relaxation} of $E$. More specifically, we refer to $E$ as a mixed-integer linear programming (MILP) (resp. mixed-integer convex programming (MICP)) relaxation if $R$ is a polyhedron (resp. convex set).  Moreover, we say that $E$ is an \textit{MIP formulation} of $S$ if $E$ models $S$ exactly, that is, $S = \proj_{x}(E)$. We will discuss the \textit{quality} of a given MIP relaxation $E$ of $S$ in terms of the size $(|y|, |\delta|)$ as well as the strength of the associated continuous relaxation $R$. We say that an MIP relaxation $E$ is \textit{ideal} if $\proj_\delta\bigl(\vertex(R)\bigr) \subseteq \{0,1\}^{n_3}$~\cite{vielma2015mixed}. An ideal formulation is desirable since an extreme optimal solution to the continuous relaxation is integer and, therefore, optimal for the MIP formulation. Moreover, branching on any binary variables restricts $E$ to a face, and, thus, ideality is retained upon branching.

To construct MIP relaxations for the hypograph of $\phi \mcirc f$, we use a vector $a$ to partition the hypercube $[f^L,f^U]$, which contains the range of $f(\cdot)$ over  $X$. 
%This vector captures the discretization points. 
\revision{Formally, for $i  \in [d]$, let $n_i+1$ be the number of points used to partition the interval $[f_i^L,f_i^U]$.  Throughout this paper, we will assume without loss of generality and for notational simplicity that $n_1 = \cdots = n_d = n$, where $n$ is a positive integer. To present the partition, we consider a vector $a := (a_1, \ldots, a_d) \in \R^{d \times (n+1)}$ so that, for  $i \in [ d ]$, $f^L_i = a_{i0} < \cdots < a_{in} = f^U_i$. The grid partition points, $\prod_{i=1}^d\{a_{i0}, \ldots, a_{in}\}$, play two roles in our construction. First, each point $a_{ij}$  will be related to an underestimator of $f_i(\cdot)$ to capture inner-function structure. Second, these values will be used to discretize the outer-function. 
%First, as detailed in Section~\ref{section:MIP-conv}, a subset of them serve as the discretization points. Second, we use them to capture inner-function structure.
However, to control the number of binary variables in the relaxation, we will discretize the $i^{\text{th}}$ coordinate only at a subsequence $a_{i\tau(i,0)}, \ldots, a_{i\tau(i,l_i)}$ of $a_{i0}, \ldots, a_{in}$, where $0=\tau(i,0)< \ldots < \tau(i,l_i) = n$.} 
% This is done to control the number of binary variables in the relaxation. 
 As a result, we obtain a subdivision $\mathcal{H}$ of $[f^L,f^U]$,
\[
\mathcal{H}: = \Bigl\{ \prod_{i = 1}^d\bigl[a_{i\tau(i,t_i-1)},a_{i\tau(i,t_i)}\bigr] \Bigm| t_i \in [l_i] \text{ for each } i \in [d] \Bigr\}. 
\]
Henceforth, this pair $(a, \mathcal{H})$ will be referred to as a \textit{discretization scheme} for a composite function $\phi \mcirc f$.

Given such a discretization $(a,\mathcal{H})$, the typical strategy of constructing an MIP relaxation for the hypograph of $\phi \mcirc f$ is to outer-approximate the graph of the inner-function $f(\cdot)$ with a polyhedron, and  \revision{construct an MIP formulation of the following disjunctive constraints}
\begin{equation}\label{eq:MIP-standard}
	 (f,\phi) \in \bigcup_{H \in \mathcal{H}} \conv\bigl( \hypo(\phi|_H)\bigr),
\end{equation}
where $\phi|_H(\cdot)$ denotes the restriction of $\phi(\cdot)$ on $H$. To derive an MIP formulation for~(\ref{eq:MIP-standard}), we can construct $\conv(\phi|_H)$ for each $H \in \mathcal{H}$ and then use disjunctive programming~\cite{balas1998disjunctive}. Using this approach, the standard formulation for~(\ref{eq:MIP-standard}) from~\cite{balas1998disjunctive,jeroslow1984modelling} requires a binary variable $\delta^H$ and a copy of $(f^H, \phi^H)$ for $\conv(\phi^H)$, thus introducing $|\mathcal{H}|$  auxiliary binary variables and $|\mathcal{H}|(d+1)$ auxiliary continuous variables. Although this formulation of~(\ref{eq:MIP-standard}) is ideal, the formulation and its continuous relaxation are typically intractable since $|\mathcal{H}|$ is exponential in $d$. 

In Section~\ref{section:MIP-conv}, we provide a new MIP relaxation framework for the hypograph of $\phi \mcirc f$ which seamlessly integrates the incremental model~\cite{dantzig1960significance}. Then in Section~\ref{section:MIP-relaxation}, we show in \revision{Theorem~\ref{them:DCR}} that this approach allows us to leverage the recently introduced composite relaxation framework~\cite{he2021new}. A key feature is that our scheme uses existing envelope results alongside ideal formulations of a combinatorial set to derive ideal MIP formulation for \eqref{eq:MIP-standard}. Additionally, due to the relationship with the composite relaxation framework, we are able to exploit inner-function structure. This allows us to derive  formulations that are simultaneously tighter and more compact than the exponentially sized disjunctive programming based counterparts. As an example, \revision{in Corollary~\ref{cor:MIP-phi-supermodular}}, we show that when the envelope of a specific lifting of the outer-function $\phi(\cdot)$ can be separated easily, our MIP formulation for~(\ref{eq:MIP-standard}) requires $d(n+1)$ auxiliary continuous variables and the continuous relaxation is tractable.
% \revision{\sout{The number of binary variables depends on the formulation for the combinatorial set and could be as few as $\sum_{i=1}^d\log_2(l_i-1)$}}. 
To the best of our knowledge, for this setting, none of the previous ideal MIP formulations had polynomially many continuous variables or a tractable continuous relaxation.  Section~\ref{section:comparison} will provide geometric insights into the quality of relaxations, paving the way for tightening our MIP relaxations. Section~\ref{section:computation} compares our relaxations with those used in current state-of-the-art solvers and shows that the relaxations are able to close 70\% of gap on our test instances.
%Section~\ref{section:MIP-log} will be devoted to reduce the number of binary variables of our MIP relaxations to $\mathcal{O}(\sum_{i=1}^d \log l_i)$  without increasing the size of continuous variables.

\section{Ideal MIP formulations via convexification}\label{section:MIP-conv}
For each $i \in [ d ]$, let $\Delta_i: = \{ z_i \in \R^{n+1} \mid 1 = z_{i0} \geq z_{i1} \geq \cdots \geq z_{in} \geq 0\}$, and let $\Delta$ be the \textit{simplotope} defined as $\Delta := \Delta_1 \times \cdots \times \Delta_d$. Then, a model of selecting subcubes from $\mathcal{H}$ is given as follows:
\begin{subequations}\label{eq:Inc}
\begin{align}
		& 
	z_i \in \Delta_i \quad \delta_{it} \in \{0,1\}\quad z_{i\tau(i,t)} \geq \delta_{it} \geq z_{i\tau(i,t)+1} &&  \text{for }i \in [ d ]\text{ and } t \in [ l_i-1 ],
 \label{eq:Inc-1}  \\
		&f_i =  a_{i0}z_{i0} + \sum_{j = 1}^{n}(a_{ij} - a_{ij-1})z_{ij} =: F_i(z_i) &&  \text{for }i \in [  d ]  \label{eq:Inc-2}. 
\end{align}	
\end{subequations}
\revision{For each $i \in [d]$ and $j \in [n]$, the variable $z_{ij}$ can be interpreted as the portion of interval length $a_{ij} - a_{ij-1}$ that is selected to define variable $f_i$. Moreover, each binary variable $\delta_{it}$ and its relation with two adjacent variables $z_{i\tau(i,t)}$ and $z_{i\tau(i,t)+1}$ forces $f_i$ to be defined in an incremental fashion.} \revision{With each hypercube $H = \prod_{i=1}^d[a_{i\tau(i,t_i-1)},a_{i\tau(i,t_i)}] \in \mathcal{H}$, we associate a face of $\Delta$, that is, 
\[
\Delta_H:=\bigl\{z \in \Delta \bigm| z_{ij} = 1 \text{ for } j \leq \tau(i,t_i-1) \text{ and } z_{ij} = 0 \text{ for } j> \tau(i,t_i) \; \forall i \in [d] \bigr\}.
\]}
In order to analyze our MIP relaxations, we will use the following lemma that relates the facial structure of $\Delta$ to the collection of hypercubes $\mathcal{H}$ introduced above. \revision{We omit the proof of the following result because it follows by direct computation.}
\begin{lemma}\label{lemma:Inc}
An MILP formulation for $\{ \Delta_H\}_{H \in \mathcal{H}}$ is given by~(\ref{eq:Inc-1}), \revision{that is $z \in \{ \Delta_H\}_{H \in \mathcal{H}}$ if and only if there exists $\delta$ such that $(z,\delta)$ satisfies~(\ref{eq:Inc-1}).}  Let $F(z):=\bigl(F_1(z_1), \ldots, F_d(z_d)\bigr)$ for every $z \in \Delta$. Then, for each  $H \in \mathcal{H}$, $F(\Delta_H) = H$.\end{lemma}
% \begin{compositeproof}
% The proof is omitted.\Halmos
% \end{compositeproof}
\begin{remark}
\revision{If all the points are used for discretization then the index $\tau(i,t)$ in~\eqref{eq:Inc-1} simplifies to $t$. If $d$ is further assumed to be $1$ then the formulation~(\ref{eq:Inc}) is the  \textit{incremental formulation}~\cite{dantzig1960significance} of selecting intervals from $\bigl\{[a_{1j-1}, a_{1j}]\bigr\}_{j = 1}^n$. Henceforth,~(\ref{eq:Inc}) will be referred to as the incremental formulation that selects subcubes from $\mathcal{H}$.}\QED
\end{remark}
%We remark that when $d=1$ and $\{\tau(i,t)\}_{t = 0}^{l_i} =  \{0\} \cup [ n ]$, the formulation~(\ref{eq:Inc}) is the  \textit{incremental formulation}~\cite{dantzig1960significance} of selecting intervals from $\bigl\{[a_{1j-1}, a_{1j}]\bigr\}_{j = 1}^n$.
%\[
%f_1 = F_1(z_1) \quad z_1 \in \Delta_1 \quad   \delta_1 \in \{0,1\}^{n-1} \quad \text{and }  z_{1j} \geq \delta_{1j} \geq z_{1j+1}\quad \text{for } j \in [  n-1 ] .
%\]
\revision{
The incremental formulation is often used to model piecewise-linear functions~\cite{vielma2011modeling,yildiz2013incremental,vielma2015mixed} and construct MIP relaxations for univariate nonlinear functions~\cite{geissler2011using}. It has been observed computationally in~\cite{vielma2011modeling,yildiz2013incremental} that incremental formulations often require fewer branch-and-bound nodes to solve the problem. Our usage of incremental formulations will be different from the literature in two ways. First, we  show in Theorem~\ref{them:MIP-phi-z} that, under certain conditions, the concave envelope  of a multivariate outer-function $\phi(\cdot)$ over $(f,z)$  yields an ideal formulation for \eqref{eq:MIP-standard}. 
%This is useful since we do not need to consider $\delta$ variables during the convexification. 
Later, in Theorem~\ref{them:DCR}, we relate the $z_{ij}$ variable in the incremental formulation to the inner-function structure for a composite function. }
%the slope of the line connecting $(a_{ij-1}, s_{ij-1})$ to $(a_{ij},s_{ij})$, where $s_{ij}$ represents a particular under-estimator of the inner-function $f_i(\cdot)$ bounded from above by $a_{ij}$.}

% is the evaluation, at any $x$, of an underestimator of inner function bounded by $a_{ij}$. 
%We use this insight in Theorem~\ref{them:DCR} to tighten the relaxation by using the inner-function structure via its underestimators.

Next, we provide an ideal formulation for~(\ref{eq:MIP-standard}). Consider a function $\phi \mcirc F: \Delta \to \R$ defined as $(\phi \mcirc F)(z) = \phi\bigl(F_1(z_1), \ldots, F_d(z_d) \bigr)$, where $F_i(\cdot)$ is introduced in the incremental formulation~(\ref{eq:Inc-2}). 
%The main idea of our derivation is to utilize the concave envelope of $\phi \mcirc F$ over $\Delta$. 
We will show that
\begin{equation}\label{eq:MIP-z}
	\Bigl\{(f,\phi,z,\delta) \Bigm|  \phi \leq \conc_\Delta (\phi \mcirc F)(z)  ,\ (f,z,\delta) \in(\ref{eq:Inc}) \Bigr\}
\end{equation}
is an MIP formulation for~(\ref{eq:MIP-standard}). Moreover, this formulation is ideal if $\phi \mcirc F$ satisfies a certain convex extension property defined as follows. 
\begin{definition}[\cite{tawarmalani2002convex}]
	A function $g: S \to \R$, where $S$ is a polytope, is said to be concave-extendable (resp. convex-extendable) from $Y \subseteq S$ if the concave (resp. convex) envelope $g$ and $g|_Y$ over $S$ are identical, where $g|_Y$ is the restriction of $g$ to $Y$. To assess concave-extendability (resp. convex-extendability), $g|_Y(x)$ is defined as $g(x)$ for $x\in Y$ and $-\infty$ (resp. $\infty$) otherwise.\QED
\end{definition} 
Our proof will need the following lemma, regarding ideality of MIP formulations.  
\begin{lemma}\label{lemma:ideal}
Let $S$ be a polyhedron in $\R^{n+m}$ such that for $1\le j,k\le m$ 
\[
(x,y) \in  \vertex(S) \cup \vertex\bigl(S \cap\{(x,y) \mid y_{j} = y_k\} \bigr) \implies y \in \mathbb{Z}^m.
\]
For $\hat{S} = \{(x,y,\delta)\in S \times \R\mid y_{j} \leq \delta\leq y_k \} \neq \emptyset$, we have $  \proj_{(y,\delta)}\bigl(\vertex(\hat{S})\bigr) \subseteq \mathbb{Z}^{m+1}$.
\end{lemma}

\begin{compositeproof}
Consider an extreme point $(\hat{x},\hat{y}, \hat{\delta})$ of $\hat{S}$. It is easy to see that either $\hat{\delta} = \hat{y}_j$ or $\hat{\delta} = \hat{y}_k$. Therefore, it suffices to show that $\hat{y} \in \mathbb{Z}^m$. Without loss of generality, assume that $\hat{\delta} =\hat{y}_j$. There are two cases to consider, either  $\hat{y}_j < \hat{y}_k$ or $\hat{y}_j = \hat{y}_k$. Let $\hat{y}_j < \hat{y}_k$. Then, since $(\hat{x},\hat{y},\hat{\delta})$ is an extreme point of $\hat{S}$ and, therefore of $\hat{S}\cap\{(x,y,\delta\mid \delta=y_j\}$, it follows that $(\hat{x}, \hat{y})$ is an extreme point of $S$ since the latter set is obtained by an injective affine transformation. It follows from the hypothesis that $\hat{y} \in \mathbb{Z}^m$. %If not, $(\hat{x}, \hat{y})$ is expressible as a convex combination of  $( x',y')$ and $(x'',y'')$ of $S$. Then, for sufficient small $\epsilon > 0$, $(\hat{x},\hat{y}) + \epsilon(x'-x'', y' - y'')$ and $(\hat{x},\hat{y}) - \epsilon(x'-x'', y' - y'')$ are in $\proj_{(x,y)}(\hat{S})$. 
If $\hat{y}_j = \hat{y}_k$, we claim that $(\hat{x}, \hat{y})$ is an extreme point of $S \cap \{ (x,y) \mid y_j = y_k \}$, and, thus, by the hypothesis $\hat{y} \in \mathbb{Z}^m$. If not, $(\hat{x},\hat{y})$ is expressible as a convex combination of two distinct points $( x',y')$ and $(x'',y'')$ of $S \cap\{ (x,y) \mid y_j = y_k \}$. Then, $(\hat{x}, \hat{y}, \hat{\delta})$ is a convex combination of $( x',y',y'_j)$ and $( x'',y'',y''_j)$ of $\hat{S}$, yielding a contradiction. \Halmos
\end{compositeproof}

%\begin{compositeproof}
%Let $(x,y,\delta)$ be an extreme point of $\hat{R}$. First, we observe that $\delta=y_j$ or $\delta = y_k$. Suppose that this is not the case. Then, there exists $\epsilon > 0$ such that points $(x,y,\delta -\epsilon)$ and $(x,y,\delta +\epsilon)$ belong to $\hat{R}$, a contradiction. Without loss of generality, we assume that $\delta = y_j$. Next, we show that $y \in \{0,1\}^{n_2}$ by contradiction. This implies that ${0,1}^{n_2 +1}$. Suppose that $y \notin \{0,1\}^{n_2}$. Then, $(x,y)$ is not an extreme point of $R$ since $\proj_y\bigl(\vertex(R) \bigr)$ is assumed to be in $ \{0,1\}^{n_2}$,. It follows readily that there exist convex multipliers $\{\gamma^t\}_{t \in T}$ and $\{x^t,y^t\}_{t \in T}$ such that $(x,y,\delta) = \sum_{t \in T} \gamma_t(x^t,y^t,y^t_j)$. This yields a contradiction as $\{(x^t,y^t,y^t_j)\}_{t \in T} \subseteq \hat{R}$, \Halmos
%\end{compositeproof}

\begin{theorem}\label{them:MIP-phi-z}
An MIP formulation for~(\ref{eq:MIP-standard}) is given by~(\ref{eq:MIP-z}). If $\phi \mcirc F: \Delta \to \R $ is concave-extendable from $\vertex(\Delta)$ then~(\ref{eq:MIP-z}) is an ideal MILP formulation.  
\end{theorem}
\begin{compositeproof}
First, we show that~(\ref{eq:MIP-z}) is an MIP formulation of~(\ref{eq:MIP-standard}).  For $C \subseteq \Delta$, let $E_{C}:= \bigl\{ (f,\phi,z) \bigm| z \in C,\ f = F(z),\ \phi \leq \phi(f) \bigr\}$. By linearity of $F(\cdot)$, it follows that $\conv(E_\Delta) = \bigl\{(f,\phi,z)\bigm| z \in \Delta,\ f = F(z),\ \phi \leq \conc_{\Delta}(\phi \mcirc F )(z)\bigr\}$.
Then, by the first statement in Lemma~\ref{lemma:Inc}, projecting out the binary variables from~(\ref{eq:MIP-z}) yields $\bigl\{ \conv(E_\Delta) \cap (\R^{d+1} \times \Delta_H)\bigr\}_{H \in \mathcal{H} }$. Therefore, the proof is complete since projection commutes with set union and for $H \in \mathcal{H}$
\[
\begin{aligned}
\proj_{(f,\phi)}\bigl( \conv(E_\Delta) \cap (\R^{d+1} \times \Delta_H) \bigr)& = \proj_{(f,\phi)}\bigl(\conv(E_{\Delta_H})\bigr) \\
& = \conv\bigl(\proj_{(f,\phi)}(E_{\Delta_H})\bigr) = \conv\bigl(\hypo(\phi|_H)\bigr),
\end{aligned}
\]
where the first equality holds since $\Delta_H$ is a face of $\Delta$, the second equality holds as the projection commutes with convexification, the last equality holds since, for $C \subseteq \Delta$, $\proj_{(f,\phi)}(E_C) = \hypo(\phi|_{F(C)})$, and, by Lemma~\ref{lemma:Inc}, $F(\Delta_{H}) =H$. 

Next, we prove the second statement. Suppose that $\phi \mcirc F$ is concave-extendable from $\vertex(\Delta)$. It follows readily that $\conc_\Delta(\phi \mcirc F)(\cdot)$ is a polyhedral function, and, thus, the MIP formulation~(\ref{eq:MIP-z}) is linear. Now, we apply Lemma~\ref{lemma:ideal} recursively to show that the formulation is ideal. For $t = (t_1, \ldots, t_d) \in \prod_{i=1}^d\{0\} \cup [ l_i-1 ]$, let
\[
\begin{aligned}
R^t: = \Bigl\{(f, \phi,z, \delta) \Bigm|&  (f, \phi, z) \in \conv(E_\Delta),\ \delta \in [0,1]^{\sum_{i=1}^d(l_i-1)},  \\
 & \qquad  z_{i\tau(i,j_i)+1} \leq \delta_{ij_i} \leq z_{i\tau(i,j_i)} \text{ for }  i \in [  d ] \text{ and } j_i \in  [ t_i ] \Bigr\},
\end{aligned} 
\]
and let $y^t: = \bigl((\delta_{11}, \ldots, \delta_{1t_1} ), \ldots, (\delta_{d1}, \ldots, \delta_{dt_d} )\bigr)$. Since  $\phi \mcirc F$ is concave-extendable from $\vertex(\Delta)$, $\proj_{z}\bigl(\vertex(R^0)\bigr) = \vertex(\Delta)$, a set consisting of binary points. Let $\iota \in [ d ]$ and $t' \in \prod_{i=1}^d \{0\} \cup [ l_i-1 ]$ so that $t'_\iota < l_\iota-1$, and assume that the points in $ \proj_{(z,y^{t'})}\bigl(\vertex(R^{t'})\bigr)$ are binary. Since $z_{\iota\tau(\iota,t'_\iota+1)}= z_{\iota\tau(\iota,t'_\iota+1)+1}$ defines a face $F$ of $R^{t'}$, $\vertex(F) \subseteq \vertex(R^{t'})$, and thus, by hypothesis, $\proj_{(z,y^{t'})}\bigl( \vertex(F) \bigr)$ is binary. By Lemma~\ref{lemma:ideal}, $\proj_{(z,y^{t''})}\bigl(\vertex(R^{t''})\bigr)$ is binary, where $t''_i = t'_i$ for $i \neq \iota$ and $t''_\iota = t'_\iota+1$.  \Halmos
\end{compositeproof}
\begin{remark}
It follows from the proof of Theorem~\ref{them:MIP-phi-z} that a formulation for 
$\cup_{H \in \mathcal{H}} \conv\bigl(\graph(\phi|_H)\bigr)$
can be obtained by convexifying the graph of $\phi \mcirc F$ over $\Delta$, \revision{ 
\[
\bigl\{(f,\phi,z,\delta) \bigm|  (\phi,z) \in  \conv\bigl(\graph(\phi \mcirc F)\bigr),\ (f,z,\delta) \in(\ref{eq:Inc}) \bigr\},
\]
which is ideal if $\phi \mcirc F$ is convex- and concave-extendable from $\vertex(\Delta)$.} \QED
\end{remark}
\revision{In Proposition~\ref{prop:separateIntegral}, we show that replacing~\eqref{eq:Inc-1} with any formulation of $\{\Delta_H\}_{H \in \mathcal{H}}$ yields an MIP formulation for~\eqref{eq:MIP-standard} and discuss its ideality.} 
% we generalize Theorem~\ref{them:MIP-phi-z} to allow the usage of any formulation of $\{\Delta_H\}_{H \in \mathcal{H}}$. 
The next example shows that the extension condition can not be relaxed.
\begin{example}~\label{ex:ideal-extension}
Consider a function $\phi(\cdot)$ defined as $\phi(f) = \sqrt{f}$ for $0\leq f \leq 1$ and $\phi(f)= f$ for $1< f\leq 2$. Choosing $(a_0,a_1,a_2) = (0,1,2)$ as discretization points, formulation~(\ref{eq:MIP-z}) yields
\begin{equation}\label{eq:MIP-z-ex}
\phi \leq \conc_\Delta(\psi)(z) \quad 1 = z_0 \geq z_1 \geq \delta \geq z_2 \geq 0 \quad \delta \in \{0,1\} \quad \text{and } f = z_1 + z_2,
\end{equation}
where $\Delta:=\{z \mid 1= z_0 \geq z_1 \geq z_2 \geq 0\}$ and $\psi(z) = \phi(z_1 + z_2)$. Notice that the composite function $\psi(\cdot)$ does not satisfy the extension condition since $\psi(\cdot)$ is concave over the edge $\{z \mid 1= z_0 \geq z_1 \geq z_2 = 0\}$. Now, \revision{let $f = \frac{1}{2}$, $\phi = \frac{1}{2}$, $z = (z_0,z_1,z_2) =  (1,\frac{1}{2},0)$ and $\delta = \frac{1}{2}$}. It can be verified that $(f,\phi,z,\delta)$ is extremal, proving that~(\ref{eq:MIP-z-ex}) is not ideal.\QED
\end{example}

When the extension condition does not hold, it is possible to obtain an ideal formulation by constructing the concave envelope in the extended space of $(z,\delta)$ variables. More specifically, let $\hat{\Delta}: = \hat{\Delta}_1 \times \cdots \times \hat{\Delta}_d$, where $\hat{\Delta}_i$ is defined as:
\begin{equation}\label{eq:extended-Z}
\hat{\Delta}_i: = \bigl\{(z_i,\delta_i) \bigm| z_i \in \Delta_i,\ z_{i\tau(i,t)} \geq \delta_{it} \geq z_{i\tau(i,t)+1} \text{ for } t \in [  l_i-1 ]  \bigr\}.
\end{equation}
Moreover, let $\widehat{\phi F}: \hat{\Delta} \to \R$ be an extension of $\phi\circ F$ defined as $\widehat{\phi F}(z,\delta) = (\phi\circ F)(z)$ for $(z,\delta) \in \hat{\Delta} \cap \bigl(\R^{d \times (n+1)} \times \{0,1\}^{\sum_{i=1}^dl_i-1}\bigr)$ and $\widehat{\phi F}(z,\delta) = \infty$ otherwise.  
\begin{proposition}\label{prop:zdeltaenv}
An ideal MIP formulation for~(\ref{eq:MIP-standard}) is given by
\begin{equation}\label{eq:zdeltaformulation}
\Bigl\{(f,\phi,z,\delta) \Bigm|  \phi \leq \conc_{\hat{\Delta}} (\widehat{\phi F})(z,\delta)  ,\ (f,z,\delta) \in(\ref{eq:Inc}) \Bigr\}. 
\end{equation}
\end{proposition}
% \begin{compositeproof}
% \revision{Provide a proof here}. \Halmos 
% \end{compositeproof}
\begin{compositeproof}
\revision{The validity of~\eqref{eq:zdeltaformulation} follows from that of~\eqref{eq:MIP-z}. More specifically, for each hypercube $H = \prod_{i=1}^d[a_{i\tau(i,t_i-1)},a_{i\tau(i,t_i)}] \in \mathcal{H}$, let $\hat{\Delta}_H$ be the face of $\hat{\Delta}$ such that, for $i \in \{1, \ldots,d \}$,  $\delta_{ij} = 1$ for $j < t_i$ and $\delta_{ij} = 0$ for $j \geq t_i$ and  $z_{ij} = 1$ for $j \leq \tau(i,t_i-1)$  and $z_{ij} = 0$ for $j> \tau(i,t_i)$. Then, $(z,\delta)$ satisfies~(\ref{eq:Inc-1}) if and only if there exists a hypercube $H \in \mathcal{H}$ so that $(z,\delta) \in \hat{\Delta}_{H}$ and
\[
\conc_{\hat{\Delta}}( \widehat{\phi F})(z,\delta) = \conc_{\hat{\Delta}_H}(  \widehat{\phi F})(z,\delta) = \conc_{\Delta_H}( \phi \mcirc F)(z) = \conc_{\Delta}( \phi \mcirc F)(z),
\]
where the first (resp. third) equalities hold because $\hat{\Delta}_H$ (resp. $\Delta_H$) is a face of $\hat{\Delta}$ (resp. $\Delta$), and the second equality follows since, for fixed $\delta$, $ \widehat{\phi F}(z) = (\phi \mcirc F)(z)$.}

\revision{
To show the ideality of \eqref{eq:zdeltaformulation}, consider a vertex $(f,\phi,z,\delta)$ of its LP relaxation, and assume that $\delta$ is not binary. Without loss of generality, assume $\phi =\conc_{\hat{\Delta}} (\widehat{\phi F})(z,\delta)$.  Notice that $(z,\delta) \notin \hat{\Delta}':=\cup_{H \in \mathcal{H}}\hat{\Delta}_H$ and $\hat{\Delta} = \conv(\hat{\Delta}')$. 
Thus, $(z,\delta, \phi)$ is a convex combination of  points $\bigl(z^k,\delta^k, \phi^k\bigr)$, where $(z^k,\delta^k) \in \hat{\Delta}'$ and $\phi^k = \conc_{\hat{\Delta}}(\widehat{\phi F})(z^k,\delta^k)$. Therefore, $(f,\phi,z,\delta)$ can be expressed as a convex combination of distinct points $\bigl(F(z^k),z^k, \delta^k,\phi^k \bigr)$ feasible to \eqref{eq:zdeltaformulation}, a contradiction to extremality of $(f,\phi,z,\delta)$.} \Halmos
\end{compositeproof}
\revision{As a result, Proposition~\ref{prop:zdeltaenv} yields an ideal formulation for the previous example: 
\[
\phi \leq \sqrt{(1-\delta)(z_1 - \delta)} +z_2,\ 1 = z_0 \geq z_1 \geq \delta \geq z_2 \geq 0,\ \delta \in \{0,1\}, \ f = z_1 +z_2.
%,
\]
%where the envelope description is obtained using disjunctive programming~\cite{balas1998disjunctive}.
It can be verified that $(\phi,z,\delta,f) = (\frac{1}{\sqrt{2}},1,\frac{1}{2},0, \frac{1}{2}, \frac{1}{2})$ is not feasible to the above constraints, certifying that formulation~(\ref{eq:MIP-z-ex}) is not ideal. For the remainder of this paper, we focus on functions $\phi \mcirc F$ that are concave-extendable, and therefore, it will suffice to construct the concave envelope description in $z$-space. }

%For the remainder of this paper, we focus on functions $(\phi \mcirc F)(\cdot)$ that are concave-extendable, and thus, it will suffice to derive the envelope description in $z$-space.
% \begin{example}
% Consider the setting in Example~\ref{ex:ideal-extension} again. Here, Proposition~\ref{prop:zdeltaenv} yields an ideal formulation given as follows: 
% \[
% \begin{aligned}
% &\phi \leq \sqrt{(1-\delta)(z_1 - \delta)} +z_2 \quad 1 \geq z_1 \geq \delta \geq z_2 \geq 0 \quad  \delta \in \{0,1\}  \text{ and }  f = z_1 +z_2
% \end{aligned}
% \]
% where the envelope description over the space of $(z,\delta)$ variables is obtained using disjunctive programming~\cite{balas1998disjunctive}. It can be verified that $(f,\phi,z,\delta) = (\frac{1}{2},\frac{1}{\sqrt{2}},1,\frac{1}{2},0, \frac{1}{2})$, the point used in Example~\ref{ex:ideal-extension}, is not feasible to the above constraints.    \QED \end{example}
%\begin{proposition}[ideal formulations for piecewise-concave functions]
%\end{proposition}

\section{MIP relaxations for composite functions}\label{section:MIP-relaxation}

\revision{In this section, we leverage Theorem 3.4 and the composite relaxation framework \cite{he2021new} to develop new MIP discretization-based relaxations for mixed-integer nonlinear programs. After providing a brief overview of the relaxation technique, we propose the relaxation in Theorem~\ref{them:DCR}. Then, we specialize the results to specific classes of outer functions, for example, we treat bilinear outer function in the space of original variables in Corollary~\ref{eq:bilinear-stair-2} and multilinear outer-functions in Corollary~\ref{cor:mul}.}

%We then provide insights into why composite relaxations improve the quality of relaxations and provide a further tightening that exploits local bounds on estimators for inner functions. 

\subsection{MIP relaxations via exploiting inner-function structure}
%The grid partition points, $\prod_{i=1}^d\{a_{i0}, \ldots, a_{in}\}$, play a dual role in our MIP relaxations.  First, as detailed in Section~\ref{section:MIP-conv}, a subset of them serve as the discretization points. Second, we use them to capture inner-function structure. Towards this end, 
%In this subsection, we develop new insights into the recently developed composite relaxations (CR) \cite{he2021new}. After a brief review, we show that this allows us to substantially improve tightness of the disjunctive programming based MIP as well as CR relaxations.
We start with reviewing the recently developed composite relaxations (CR) \cite{he2021new}. For $i \in [ d ]$, assume that there exists a vector of functions $u_i: X \to \R^{n+1}$ such that for every $x \in X$ 
\begin{equation}\label{eq:ordered-oa}
u_{i0}(x) = a_{i0} \quad u_{in}(x) = f_i(x)\quad \text{and }	u_{ij}(x) \leq \min \bigl\{f_i(x), a_{ij} \bigr\}  \text{ for } j \in [ n-1 ].
\end{equation}
In other words, function $u_{ij}(\cdot)$ is an underestimator of the inner-function $f_i(\cdot)$ which is bounded from above by the point $a_{ij}$. Then, the ordering relationship in~(\ref{eq:ordered-oa}) is linearized into a polytope, $P:=\prod_{i=1}^d P_i$,  where 
\begin{equation}\label{eq:P-3}
P_i = \left\{ u_i \in \R^{n+1} \,\middle|\,
  \begin{aligned}
&a_{i0} \leq u_{i j} \leq  \min\{a_{i j}, u_{in}\} \text{ for } j \in [ n-1] \\ &u_{i 0} = a_{i 0},\; a_{i 0} \leq u_{i n} \leq a_{i n}
  \end{aligned}
  \right\}.
\end{equation}
Let $\ephi: P \to \R$ be an extension of the outer-function $\phi:[f^L, f^U] \to \R$ defined as $\ephi(u) = \phi(u_{1n}, \ldots, u_{dn})$ for every $u \in P$. Observe that, in $P$, there is a $u$ variable for each underestimator while $\bar{\phi}(\cdot)$ depends only on a few coordinates of $u$, those \revision{representing} inner-functions $f(\cdot)$. Then, the \textit{composite relaxation} for  $\phi \mcirc f$ is obtained by convexifying the hypograph of $\ephi(\cdot)$ over the polytope $P$, that is 
\begin{equation}\label{eq:CR}
\hypo(\phi \mcirc f) \subseteq \bigl\{(x,\phi) \bigm| \phi \leq \conc_P(\ephi)(u),\ u(x)\leq u,\ (x,u_{\cdot n}) \in W \bigr\},
\end{equation}
where $u_{\cdot n}=(u_{1n}, \ldots, u_{dn})$  and $ \graph(f) \subseteq W$.
It was shown in \cite{he2021new} that \revision{the} convex hull over $P$ is expressible in extended space with a few extra variables if the convex hull over a subset $Q$ of $P$ is known. Here, $Q:= \prod_{i=1}^dQ_i$ and $Q_i$ is a simplex in $\R^{n+1}$ with the following extreme points,
\begin{equation}\label{eq:Q-V}
v_{ij} = (a_{i0}, \ldots, a_{ij-1}, a_{ij}, \ldots, a_{ij} ) \qquad \text{for } j \in  0 \cup [ n ] . 
\end{equation} 
The extended formulation  of the hypograph of $\conc_P(\ephi)(\cdot)$ is given as follows:
\begin{equation}\label{eq:extended-P}
\bigl\{(u, \phi)\bigm| \phi \leq \conc_Q(\ephi)(s) ,\  u \in P,\ u\leq s,\ u_{\cdot n} = s_{\cdot n} \bigr\}.
\end{equation}
\revision{Recall that the variable $u$ represents underestimators $u(\cdot)$. Here, each variable $s_{ij}$ can be interpreted as a representation of an underestimator of $f_i(\cdot)$ that is bounded from above by $a_{ij}$ but is tighter than $u_{ij}(\cdot)$. }
%\begin{equation*}\label{eq:extended-P}
%\bigl\{(u, \phi)\bigm| \phi \leq \conc_Q(\ephi)(s) ,\  u \in P,\ u\leq s,\ u_{\cdot n} = s_{\cdot n} \bigr\}.
%\end{equation*}
%
%\begin{proposition}[\cite{he2021new}]\label{prop:CR}
%Consider a vector of composite functions $\theta \mcirc f: X \to \R^\kappa$ defined as $(\theta \mcirc f)(x): = \bigl((\theta_1 \mcirc f)(x), \ldots, (\theta_\kappa \mcirc f) (x)  \bigr)$. Let $\bigl(u(x),a\bigr)$ be a pair satisfying~(\ref{eq:ordered-oa}), where $a:=(a_1, \ldots, a_d)$ so that $a_{i0} < \cdots < a_{in}$, and define $\Theta^P: = \bigl\{(u,\theta) \bigm| \theta = \theta(u_{1n}, \ldots, u_{dn}),\ u \in P \bigr\}$. Then, we obtain 
%\[
%\graph(\phi \mcirc f) \subseteq \Bigl\{ (x, \phi) \Bigm| (u,\theta) \in \conv\bigl(\Theta^P\bigr) ,\ u(x) \leq u,\  (x,u_{\cdot n}) \in W  \Bigr\},
%\]
%where $u_{\cdot n}=(u_{1n}, \ldots, u_{dn})$  and $W$ outer-approximates  $\{(x, u_{\cdot n}) \bigm| u_{\cdot n} = f(x), x\in X\}$. Moreover, for a convex relaxation $R$ of $\Theta^Q: = \bigl\{(s, \theta) \bigm| \theta = \theta(s_{1n}, \ldots, s_{dn}), s\in Q \bigr\}$, we obtain  
%\begin{equation}\label{eq:extended-P}
%\conv\bigl(\Theta^P\bigr) \subseteq  \bigl\{(u,\theta)\bigm| (s, \theta) \in R,\  u \in P,\ u\leq s,\ u_{\cdot n} = s_{\cdot n} \bigr\}, 
%\end{equation}
%where the equality holds if $R = \conv(\Theta^Q)$. \end{proposition}

\revision{Next, we develop new insights into the composite relaxations by relating the variable $z$ in the incremental formulation to the additional variable $s$ that were added in the extended formulation~\eqref{eq:extended-P}}. More specifically, Theorem~\ref{them:DCR} relates $z \in \Delta$ to variables $s \in Q$ via an invertible affine transformation. This transformation is given by $Z(s)= \bigl(Z_1(s_1), \ldots, Z_d(s_i)  \bigr)$, where $Z_i : \R^{n+1} \to \R^{n+1}$ relates $z_i$ to $s_i$ so that $Z_i(s_i) = z_i$, where
\begin{equation}\label{eq:Z_trans}
	\begin{aligned}
	 z_{i0} = 1 \qquad	\text{and} \qquad z_{ij} = \frac{s_{i j} - s_{i j-1}}{a_{ij} - a_{ij-1}} \quad \text{for } j \in [ n ].
	\end{aligned}
\end{equation}
\revision{In words, the variable $z_{ij}$ in the incremental formulation can be interpreted as the slope of the line connecting $(a_{ij-1}, s_{ij-1})$ to $(a_{ij},s_{ij})$.} The inverse of $Z$, denoted as $Z^{-1}(z):= \bigl(Z_1^{-1}(z_1), \ldots, Z_d^{-1}(z_d) \bigr)$, is given as follows:
\begin{equation}\label{eq:Z_inv}
s_{ij} = a_{i0}z_{i0}  + \sum_{k = 1}^j (a_{ ik} - a_{i k-1})z_{ik} \quad \text{for } j \in \{ 0\} \cup [n]. 
\end{equation}
%Recall that $\vertex(Q_i): = \{v_{i0}, \ldots, v_{in}\}$, where $v_{ij} = (a_{i0}, \ldots, a_{ij-1}, a_{ij} \ldots, a_{ij})$. 
\revision{The validity of our MIP relaxations is built upon the combinatorial equivalence of $\Delta$ and $Q$, or more specifically, that the faces of the two polytopes are isomorphic under the transformation $Z$.} Recall that $\Delta=\prod_{i=1}^d\Delta_i$, the vertices of $\Delta_i$ are $\zeta_{ij} = (1,\ldots,1,0,\ldots,0)$ where the first $j$ coordinates are $1$. Then, $Z_i$ maps $v_{ij}$ to $\zeta_{ij}$. Conversely, $Z_i^{-1}(\zeta_{ij}) = v_{ij}$. 
%For a hypercube $H : = \prod_{i=1}^d[a_{i\tau(i,t_i-1)}, a_{i\tau(i,t_i)}]$ in the grid partition $\mathcal{H}$, Lemma~\ref{lemma:Inc} obtained $\Delta_H$ as a lifting of $H$. Here, 
\revision{Since convexification commutes with affine maps, $\Delta_H$ is obtained as an affine transformation of $Q_H$, defined as follows}
\begin{equation*}\label{eq:Q-face}
\begin{aligned}
Q_H &:= \conv\Biggl( \prod_{i=1}^d \{v_{i \tau(i,t_i-1)}, \ldots, v_{i \tau(i,t_i)} \} \Biggr)
%= Z^{-1}\Bigl(  \conv \Bigl( \prod_{i=1}^d  \{\zeta_{i\tau(i,t_i-1)}, \ldots, \zeta_{i\tau(i,t_i)} \} \Bigr)\Bigr) 
= Z^{-1}(\Delta_{H}). 
%&= \conv\Bigl( Z^{-1}  \Bigl( \prod_{i=1}^d \{\zeta_{i\tau(i,t_{i}-1)}, \ldots, \zeta_{i\tau(i,t_i)} \} \Bigr)\Bigr)\\
%&= \conv\Bigl( \prod_{i=1}^d \{v_{i \tau(i,t_i-1)}, \ldots, v_{i \tau(i,t_i)} \} \Bigr) =:Q_{H}.
\end{aligned}
\end{equation*}
%where the first equality holds by the definition of $\Delta_H$, the second equality holds because convexification commutes with affine maps and because $Z_i^{-1}$ maps $\zeta_{ij}$ to $v_{ij}$. Conversely, $Z(Q_H) = \Delta_H$. 
%As a result, we obtain that 
%\begin{equation}\label{eq:Inc-proj-Q}
%s \in \bigcup_{H \in \mathcal{H}}Q_{H} \; \Leftrightarrow \; Z(s) \in \bigcup_{H \in \mathcal{H}}\Delta_{H} \; \Leftrightarrow	 \; \exists \delta \text{ such that } \bigl(Z(s),\delta\bigr) \in(\ref{eq:Inc-1}) .  
%\end{equation}
\begin{theorem}\label{them:DCR}
Consider a discretization scheme $(\mathcal{H},a)$, and a vector of convex function $u(\cdot)$ such that the pair $\bigl(u(\cdot),a\bigr)$ satisfies~(\ref{eq:ordered-oa}). An MICP relaxation for the hypograph of $\phi \mcirc f$ is given by:
\begin{equation*}\label{eq:DCR}
\Bigl\{ (x, \phi, s, \delta) \Bigm|
\phi \leq \conc_Q(\ephi)(s), Z(s)=z,\ \bigl(z,\delta \bigr) \in (\ref{eq:Inc-1}),\  u(x) \leq s,\ (x,s_{\cdot n}) \in W
\Bigr\},
\end{equation*}
where $W$ is a convex outer-approximation of $\bigl\{(x,s_{\cdot n}) \bigm| s_{\cdot n} = f(x),\ x\in X\bigr\}$.
\end{theorem}
%\begin{compositeproof}
%By Lemma~\ref{them:MIP-phi-z} and the affine transformation $Z$, constraints $\phi\leq \conc_Q(\phi)(s)$ and $(Z(s), \delta) \in(\ref{eq:Inc-1})$ model 
%\end{compositeproof}

\begin{compositeproof}
Let $(x, \phi) \in\hypo(\phi \mcirc f)$. We extend $(x,\phi)$ by choosing $(s,z,\delta)$ and then show that the resulting point satisfies the constraints of the proposed relaxation. For all $i$ and $j$, let $s_{ij}: = \min\bigl\{f_i(x), a_{ij} \bigr\}$ and choose  $t$, and thereby $H\in \mathcal{H}$, so that $s_{\cdot n} \in H = \prod_{i=1}^d [a_{i \tau(i,t_i-1)},a_{i \tau(i,t_i)}]$. Define $z=Z(s)$. We will show $s\in Q_H$ from which it follows that $z\in \Delta_H$. Then, we use Lemma~\ref{lemma:Inc} to pick $\delta$ so that $(z,\delta)\in \eqref{eq:Inc-1}$. Now, we show that $s\in Q_H$. Let  $j'_i$ be such that $s_{in} \in [a_{ij'_i-1},a_{ij'_i}]$. Then, $\tau(i,t_i-1)<j'_i \leq \tau(i,t_i)$ and, observe that, by our construction, $s_i = (a_{i1},\ldots,a_{ij'_i-1},s_{in},\ldots,s_{in})$. Consider vertices $v_{ij'_i-1}, v_{ij'_i} \in Q_H$ and write $s_i = \lambda  v_{ ij'_i-1} + (1-\lambda)  v_{i j'_i}$, where $\lambda = (a_{i j'_i} - s_{in})/(a_{i j'_i} - a_{i j'_i-1})$. It follows that $s \in \prod_{i} \conv\bigl(\{v_{ ij'_i-1}, v_{i j'_i}\}\bigr)\subseteq Q_H$. Finally, we show that $(x,\phi,s,\delta)$ is feasible to the proposed MICP relaxation. Since $W$ outer-approximates the graph of $f_i(x)$, $(x,s_{\cdot n}) \in W$. By (\ref{eq:ordered-oa}), $u(x) \leq s$ and we chose $\delta$ so that $(z,\delta)\in \eqref{eq:Inc-1}$. Lastly, $\phi \leq \phi(s_{\cdot n}) = \ephi(s) \leq \conc_Q(\ephi)(s)$, where the first inequality follows since $(x,\phi) \in \hypo(\phi \mcirc f)$ and $s_{\cdot n} = f(x)$, the equality is by the definition of $\ephi$, and the second inequality is because $s \in Q$.
\Halmos
\end{compositeproof}
\begin{remark}~\label{rmk:over-and-under}
%It is often the case that both over- and under-estimators are available for inner functions. 
%In fact, given an expression tree of a factorable function, Algorithm~\ref{alg:propagation} in Section~\ref{section:computation} recursively generates over- and under-estimators for each node of 
Even though we have focused on underestimators for $f$, our constructions apply to overestimators as well. 
%Algorithm~\ref{alg:propagation}, described later, generates over- and under-estimators for functions described using expression trees. This algorithm outputs
\revision{Consider a tuple} $(u_i(\cdot), a_i, \chi_i)$ for the inner function $f_i(\cdot)$, where $u_i:X \to \R^{n+1}$ is a vector of functions, $a_i$ is a vector in $\R^{n+1}$, and $\chi_i$ is the indicator vector of a subset $A_i$ of $\{0\} \cup [ n ]$ containing $\{0,n\}$ such that
	\begin{equation}\label{eq:ordered-oa-over}
\begin{aligned}
				&u_{i0}(x) = a_{i0}   \quad &&\; u_{in}(x) = f_i(x), \\
				&\text{if }  (\chi_{i})_j = 1: \quad  &&\;u_{ij}(x) \le \, \min\{f_{i}(x), a_{ij}\}  \quad \text{and} \quad u_{ij}(\cdot) \text{ is convex}, \\
				&\text{if }  (\chi_{i})_j = 0: \quad &&\max\bigl\{f_i(x),a_{ij}\bigr\} \le u_{ij}(x) \quad \text{and} \quad u_{ij}(\cdot) \text{ is concave}.\\
%				u_{in(i)}(x) = f_i(x),\qquad && a_i^L\le a_{i0} \leq \cdots \leq a_{in(i)}\le a_i^U, \\
%	 &u_{ij}(x) \leq f_{i}(x), && a_i^L \leq u_{ij}(x) \leq a_{ij} \quad \text{for } j \in A_i, \\
%  &f_{i}(x)  \leq u_{ij}(x),&& a_{ij} \leq u_{ij}(x) \leq a_i^U \quad \text{for } j \in B_i,\\
\end{aligned}
\end{equation}
A key mapping that enables the use of \eqref{eq:ordered-oa-over} is the affine transformation, $T$, introduced in~\cite{he2021new}, that maps estimators $u(\cdot)$ satisfying~(\ref{eq:ordered-oa-over}) to convex under-estimators satisfying~(\ref{eq:ordered-oa}). This transformation is defined as that $T(u)_{i j} = u_{i j}$ for $j \in A_i$ and $T(u)_{ij} = a_{i j} - u_{i j} + u_{i n}$ for $j \notin A_i$. To obtain an MICP relaxation for the hypograph of $\phi \mcirc f$, it suffices then to replace $u(\cdot)$ with $(T \mcirc u) (\cdot)$ in  Theorem~\ref{them:DCR}. \QED
\end{remark}
\begin{remark}\label{rmk:DCR-vector} Our constructions are also useful in relaxing a vector of composite functions $\theta \mcirc f$, where $\theta: \R^d \to \R^\kappa$ is a vector of functions.  Given $\bigl(u(x),a\bigr)$ that satisfies~(\ref{eq:ordered-oa}), an MICP relaxation for the graph of $\theta \mcirc f$ is given by
\begin{equation*}
\Bigl\{ (x, \theta, s, \delta) \Bigm|
(s, \theta) \in  \conv\bigl(\Theta^Q \bigr), \bigl(Z(s),\delta \bigr) \in (\ref{eq:Inc-1}),  u(x) \leq s, (x,s_{\cdot n}) \in W
\Bigr\},
\end{equation*}
where $\Theta^Q:= \{(s,\theta) \mid \theta = \theta(s_{1n}, \ldots, s_{dn}), s\in Q\}$ and $\graph(f) \subseteq W$. Since the inner functions are shared, the $s$ variables are also the same across all functions in $\theta(\cdot)$. \QED
%the set $\{(f,\theta)\mid \theta = \theta(f), f\in H\}$.
%Here, we remark that by Corollary 7 in~\cite{he2022tractable}, if, for $k \in \{1, \ldots, \kappa\}$, $\theta_k$ is supermodular over $[f^L, f^U]$ then $\conv\bigl(\hypo(\bar{\theta}|_Q) \bigr) = \cap_{k=1}^\kappa \conv\bigl(\hypo(\bar{\theta}_k|_Q) \bigr)$, where $\conv\bigl(\hypo(\bar{\theta}_k|_Q) \bigr)$ can be described using Proposition~\ref{prop:stair-Q}. 
\end{remark}

\subsection{Explicit relaxations via exploiting outer-function structure}
In this subsection, we exploit two types of outer-function structure, supermodular and multilinear, to develop explicit envelopes of $\ephi(\cdot)$ over $Q$. Then, using Theorem~\ref{them:DCR}, we write explicit MIP relaxations for the composite function $\phi \mcirc f$. For example, \revision{in Corollary~\ref{cor:MIP-phi-bilinear}}, we give a new MIP relaxation for product of non-negative functions that uses variables that scale linearly (instead of quadratically) in the number of discretizations.
\subsubsection{Supermodular outer-functions}
We start by assuming that the outer-function is supermodular. Recall that function $\eta:S \subseteq \R^n \to \R$ is said to be \textit{supermodular} if 
\[
\eta(x' \vee x'' ) + \eta(x'\wedge x'') \geq \eta(x') + \eta(x'') \qquad \text{for all }x', x'' \in S,
\]
where $x'\vee x''$ is the component-wise maximum and $x' \wedge x'' $ is the component-wise minimum of $x'$ and $x''$ and $S$ is assumed to be a lattice, that is $x'\vee x''$ and $x'\wedge x''$ belong to $S$ whenever $x'$ and $x''$ belong to $S$. Although detecting whether a function is supermodular is NP-Hard~\cite{crama1989recognition}, there are important special cases where this property can be readily detected \cite{topkis2011supermodularity}. For example, a product of nonnegative, increasing (decreasing) supermodular functions is nonnegative increasing (decreasing) and supermodular; see Corollary 2.6.3 in~\cite{topkis2011supermodularity}. 
%Also, a conic combination of supermodular functions is supermodular.
% A canonical example of a supermodular function is $\prod_{i=1}^n x_i$ over the non-negative orthant.

It is shown in Theorem 2 of~\cite{he2022tractable} that under the supermodularity condition, the concave envelope of $\ephi(\cdot)$ over $Q$ is related to the staircase triangulation of $Q$. To describe such triangulation, we associate the vertex $\text{ext}(Q,p):=(v_{1 p_1}, \ldots, v_{d p_d})$ of $Q$ with the point $p:=(p_1, \ldots, p_d)$ on a grid $\mathcal{G}:= \{ 0, \ldots, n\}^d$. A \textit{lattice path} is a sequence of points $p^{0}, \ldots, p^{r}$ in $\mathcal{G}$ such that $p^{0}= (0, \ldots, 0)$ and $p^{r}= (n, \ldots, n)$. In particular, a \textit{staircase} is a lattice path of length $dn+1$ such that for all $t \in [  dn ] $, $p^{t} - p^{t-1} = e_{k}$ where $k \in [  d ] $ and $e_{k}$ is the principal vector in $k^{\text{th}}$ direction. We refer to the movement $p^{t-1}$ to $p^{t}$ as the $t^{\text{th}}$ move. Clearly, there are exactly $n$ moves along each coordinate direction. Therefore, a staircase can be specified succinctly as $\pi:= (\pi_1, \ldots, \pi_{dn})$, where, for $t \in [ dn ] $, we let $\pi_t = i$ if $p^{t} - p^{t-1} =e_i$.  We will refer to such a vector as \textit{direction vector}, and will denote by $\Pi$ the set of all direction vectors in the grid. For $\pi \in \Pi$, let $\Upsilon_\omega$ be the simplex defined as $\Upsilon_\omega = \conv\bigl(\bigl\{ \text{ext}(Q,p_0), \ldots, \text{ext}(Q,p_{dn}) \bigr\} \bigr)$.  The set of simplices $ \{\Upsilon_\pi\}_{\pi \in \Pi}$ is called the \textit{staircase triangulation} of $Q$, that is, $Q = \cup_{\pi\in \Pi} \Upsilon_{\pi}$, and for $\pi', \pi'' \in \Pi$, $\Upsilon_{\pi'} \cap \Upsilon_{\pi''}$  is a face of both $\Upsilon_{\pi'}$ and $\Upsilon_{\pi''}$~\cite{de2010triangulations}.
 
\begin{proposition}[Theorem 2 in~\cite{he2022tractable}]\label{prop:stair-Q}
If $\ephi(\cdot)$ is concave-extendable from $\vertex(Q)$ and is supermodular over $\vertex(Q)$ then $\conc_Q(\ephi)(s) = \min_{\pi \in \Pi}\phi^\pi(s)$ for every $s \in Q$, where $\phi^\pi(\cdot)$ is obtained by affinely interpolating $\ephi(\cdot)$ over $\vertex(\Upsilon_\pi)$. \QED
\end{proposition}
\begin{remark}~\label{rmk:stair-Q}
We remark that the number of facet-defining inequalities of the hypograph of $\conc_Q(\ephi)(\cdot)$, as described in Proposition~\ref{prop:stair-Q}, is exponential in $d$ and $n$. Nevertheless, for any point $\s \in Q$, the Algorithm 1 in~\cite{he2022tractable} finds a facet-defining inequality valid for $\hypo\bigl(\conc_Q(\ephi)\bigr)$ and tight at $\s$ in $\mathcal{O}(dn \log d)$ time. \QED
\end{remark}
\begin{corollary}\label{cor:MIP-phi-supermodular}
\revision{Assume the same setup as in Theorem~\ref{them:DCR}, and assume that  $\ephi(\cdot)$ is concave-extendable from $\vertex(Q)$ and is supermodular over $\vertex(Q)$. An ideal  formulation for~\eqref{eq:MIP-standard} is  
\begin{equation*}
	\Bigl\{(\phi,f,s,\delta) \Bigm| \phi \leq \min_{\pi \in \Pi}\phi^\pi(s),\ Z(s) =z,\ f = s_{\cdot n},\ (z, \delta) \in (\ref{eq:Inc-1})\Bigr\},
\end{equation*}
which, together with $u(x) \leq s$ and $(x,s_{\cdot n}) \in W$, yields an MILP relaxation for the hypograph of $\phi \mcirc f$.}
\end{corollary}
\begin{compositeproof}
	This result follows from Theorems~\ref{them:MIP-phi-z} and~\ref{them:DCR},  Proposition~\ref{prop:stair-Q} and the invertible transformation $Z$ defined as in~\eqref{eq:Z_trans}. \Halmos
\end{compositeproof}
Corollary~\ref{cor:MIP-phi-supermodular} give new MIP relaxations for products of nonnegative increasing supermodular functions. The special case of a bilinear term is particularly interesting for three reasons. First, it is a building block in factorable relaxations~\cite{mccormick1976computability}. Second, we do not require variables to be non-negative. Third, we can derive the concave as well as convex envelope over $Q$. We will make use of this construction extensively in our computational experiments. In Corollary~\ref{cor:bilinear-stair}, we provide explicit readily implementable envelopes for the bilinear term by specializing Proposition~\ref{prop:stair-Q}. Then, Corollary~\ref{cor:MIP-phi-bilinear} gives an MIP relaxation that reduces the dimension of discretized formulations from quadratic to linear in $n$ without sacrificing quality. 
% Moreover, the bilinear term case can be used to describe envelopes of multivariate bilinear functions over $Q$. In particular, Theorem 5 in~\cite{he2022tractable} characterizes a family of bilinear functions for which both convex and concave envelopes over $Q$ can be obtained by convexifying each term in the bilinear function separately. 

%Let $U(T): u \mapsto \tilde{u}$ and $A(T): a \mapsto \tilde{a}$ so that, for $i \notin T$, $(\tilde{u}_i, \tilde{a}_i)=(u_i,a_i)$, and otherwise $\tilde{u}_{ij} =  a_{in-j} -  u_{in-j} + u_{in}$ and $\tilde{a}_{ij}= a_{in-j}$ for $j \in \{0, \ldots, n \}$.

\begin{corollary}\label{cor:bilinear-stair}
The convex envelope of $s_{1n}s_{2n}$ over $Q$ is given by 
\begin{equation}\label{eq:bilinear-stair-1}
\begin{aligned}
	\check{b}(s) := \max_{\pi \in \Pi} \biggl\{ &a_{10}a_{2n} + \sum_{t:\pi_t=1}a_{2(n-p^t_2)}\Bigl( s_{1p^{t}_1}- s_{1p^{t-1}_1} \Bigr) + \\
& \sum_{t:\pi_t=2} a_{1p^t_1}\Bigl(a_{2(n-p_2^t)} - a_{2(n-p_2^{t-1})} - s_{2(n- p^{t}_2)} + s_{2(n-p^{t-1}_2)}\Bigr)\biggr\},
\end{aligned}
\end{equation}
and the concave envelope of $s_{1n}s_{2n}$ over $Q$ is given by
\begin{equation}\label{eq:bilinear-stair-2}
 \hat{b}(s) := \min_{ \pi \in \Pi }\biggl\{ a_{10}a_{20} +  \sum_{t:\pi_t=1} a_{2p^t_2}\Bigl(s_{1p^t_1} - s_{1p^{t-1}_1}\Bigr) + \sum_{t:\pi_t=2} a_{1p^t_1}\Bigl(s_{2p^t_2} - s_{2p^{t-1}_2}\Bigr)\biggr\},
\end{equation}
where $\Pi$ is the set of all direction vectors in $\{0, \ldots, n\}^2$. 
\end{corollary} 
\begin{proof}
	See Appendix~\ref{app:bilinear-stair}.  \Halmos
\end{proof}
%We remark that discretization of bilinear terms is a common technique to relax MINLPs. By recursively applying this technique, factorable MINLPs can be approximated arbitrarily closely~\cite{misener2011apogee,misener2012global,nagarajan2019adaptive}. As such, the following result is useful in this context.
\begin{corollary}\label{cor:MIP-phi-bilinear}
Assume the same setup as in Theorem~\ref{them:DCR} and assume that the outer-function $\phi(f_1, f_2) = f_1f_2$.  Let $\check{b}$ and $\hat{b}$ be functions defined as in~(\ref{eq:bilinear-stair-1}) and~(\ref{eq:bilinear-stair-2}), respectively. An ideal  formulation for $\{ \conv(\graph(\phi|_H))\}_{H \in \mathcal{H}}$ is  
\begin{equation}~\label{eq:bilinear-MIP}
	\Bigl\{(\phi,f,s,\delta) \Bigm| \check{b}(s) \leq	\phi \leq \hat{b}(s),\ Z(s) =z,\ f = s_{\cdot n},\ (z, \delta) \in (\ref{eq:Inc-1})\Bigr\},
\end{equation}
which, together with $u(x) \leq s$ and $(x,s_{\cdot n}) \in W$, yields an MILP relaxation for the hypograph of $\phi \mcirc f$.
\end{corollary}
\begin{compositeproof}
	This result follows from Theorems~\ref{them:MIP-phi-z} and~\ref{them:DCR},  Corollary~\ref{cor:bilinear-stair} and the invertible transformation $Z$ defined as in~\eqref{eq:Z_trans}. \Halmos
\end{compositeproof}
Specifically, if there are $n-1$ discretizations along each axis, formulations from Corollary~\ref{cor:MIP-phi-bilinear} require $2(n+1)$ continuous variables and $2n$ binary variables.
% if logarithmic (resp. incremental) formulation is used. 
%These continuous variables can also be projected out using Proposition~\ref{prop:eval}, if so desired. 
This is in contrast to existing ideal formulations from~\cite{huchette2019combinatorial} which require $(n+1)^2$ continuous variables. The following example illustrates this difference. We remark that \eqref{eq:bilinear-MIP} yields tighter relaxations than those in \cite{misener2011apogee,misener2012global,nagarajan2019adaptive,huchette2019combinatorial}  because of the constraints $u(x)\le s$, and this change produces substantially tighter relaxations empirically, see Section~\ref{section:computation}.

\begin{example}~\label{ex:size-comparision}
\revision{Consider a bilinear term $f_1f_2$ over $[0,4]^2$, and a partition of the domain $\mathcal{H}:=\bigl\{ [0,3] \cup [3,4] \bigr\}^2$. In other words, we use one discretization point so that $n = 2$ and $(a_{i0}, a_{i1}, a_{i2}) = (0,3,4)$ for $i = 1, 2$. The following formulation from~\cite{huchette2019combinatorial} introduces $9$ ({\it{i.e.}}, $(n+1)^2$) continuous variables to convexify $\bigl\{(f_1, f_2, f_1f_2) \bigm| f_i \in \{a_{i0},a_{i1},a_{i2}\},\; i =1,2 \bigr\}$}:
 \[
\begin{aligned}
&(f_1,f_2,\phi) = \sum_{j=0}^2 \sum_{k = 0}^2 w_{jk}(a_{1j}, a_{2k}, a_{1j}a_{2k})\quad \sum_{j = 0}^2 \sum_{k=0}^2 w_{jk}  = 1 \quad w \geq 0, \\
	&\sum_{j=0}^2w_{j0} \leq \delta_1 \qquad \sum_{j=0}^2w_{j2} \leq 1-\delta_1 \qquad \delta_1 \in \{0,1\}, \\
	&\sum_{k=0}^2w_{0k} \leq \delta_2 \qquad \sum_{k=0}^2w_{2k} \leq 1-\delta_2 \qquad \delta_2 \in \{0,1\}.
\end{aligned}
\]
\revision{Using Corollary~\ref{cor:MIP-phi-bilinear} and  the affine map $\lambda = (T \mcirc Z)(s)$, where $Z$ is defined in~(\ref{eq:Z_trans}) and $T(z)_{ij} = z_{ij} - z_{ij+1}$, where $z_{n+1} = 0$, we introduce $6$ ({\it{i.e.}},
  $2(n+1)$) continuous variables $\lambda$  and obtain:}
\[
\begin{aligned}
&\max\left\{\!
\begin{aligned}
&12\lambda_{11}+16\lambda_{12} + 12\lambda_{21} 
%\\
%&\mskip 80mu 
+16\lambda_{22} - 16\\
& 9 \lambda_{11} + 9\lambda_{12} + 9 \lambda_{21} + 9\lambda_{22} - 9 \\
&12\lambda_{11} + 12\lambda_{12} + 9 \lambda_{21} + 12\lambda_{22} - 12 \\ 
&9\lambda_{11} +  12 \lambda_{12} + 12 \lambda_{21} + 12 \lambda_{22} - 12 \\ 
& 12 \lambda_{11} + 15\lambda_{12}  + 12 \lambda_{21} %\\
%&\mskip 80mu 
+ 15\lambda_{22} -15\\
&0 
\end{aligned}
\right\} \leq \phi \leq \min\left\{\!
\begin{aligned}
&12 \lambda_{21} + 16\lambda_{22}\\
& 3\lambda_{12} + 9 \lambda_{21} +13\lambda_{22} \\
&  4\lambda_{12} + 9 \lambda_{21} + 9\lambda_{22}\\
&12 \lambda_{11} + 16\lambda_{12}\\
& 9 \lambda_{11} +13\lambda_{12} + 3\lambda_{22} \\
&  9 \lambda_{11} + 9\lambda_{12}  4\lambda_{22}
\end{aligned}
\right\} \\
	& \lambda_i \in \Lambda_i \quad \lambda_{i0} \leq \delta_i \quad \lambda_{i2} \leq 1-\delta_i \quad \delta_i \in \{0,1\} \quad f_i = 3\lambda_{i1} + 4\lambda_{i2}  \quad \text{for } i = 1,2 .
\end{aligned}
\]
 As remarked earlier, our formulation requires significantly fewer variables than existing discretization schemes. Moreover, if $f_1=x_1^2$, we may additionally require that $3\lambda_{11}\ge \max\{2x_1-1, \frac{3}{4}x_1^2\}$ which leads to a tighter relaxation.
%Although both formulations model SOS2 constraints in the same way, we highlight that earlier discretization schemes were obtained by using disjunctive programming while our construction relies on convex hull construction over $Q$. When the convex hull is constructed using disjunctive programming or reformulation-linearization technique (RLT), our scheme recovers the earlier models. However, the new formulation can directly use convex hull formulations when they are available in a lower-dimensional setting as in Corollary~\ref{cor:MIP-phi-bilinear}. In this case, our formulation requires significantly fewer variables in comparison to disjunctive programming and RLT. In fact, our formulation can be seen as a projection of the earlier formulation. Moreover, if $f_1=x_1^2$, we have shown in Theorem~\ref{them:DCR} that we may additionally require that $3\lambda_{11}\ge \max\{2x_1-1, \frac{3}{4}x_1^2\}$ because $2x-1$ and $\frac{3}{4}x_1^2$ are convex underestimators of $x_1^2$ bounded by $3$ and the transformation $T \mcirc Z$ shows that $s_{11}=3\lambda_{11}$. One of the key advantages of Theorem~\ref{them:DCR} is that it provides a mechanism for the discretized formulation to exploit the structure of the inner functions $f_i$.
\QED
\end{example}

\subsubsection{Multilinear outer-functions} Although separating the envelope of $\ephi(\cdot)$ over $Q$ in the space of $s$ variables is NP-hard, it can often be described as the projection of a higher dimensional set with exponentially many variables. One technique that can be used in this regard is the reformulation-linearization technique ($\RLT$)~\cite{sherali1992global,sherali2013reformulation}, which we modify to convexify a multilinear function over $Q$. We begin by considering a general setting. Let $\Lambda: = \prod_{i=1}^d\Lambda_i$, where $\Lambda_i:=\bigl\{\lambda_i \in \R^{n+1} \bigm| 1- \sum_{j=0}^n \lambda_{ij}=0,\ \lambda_i \geq 0\bigr\}$. The next result shows that  the convex hull of the following multilinear set
\[
M^\Lambda: = \biggl\{(\lambda, w) \biggm| \lambda \in \Lambda,\ w_j = \prod_{i=1}^d\lambda_{ij_i} \text{ for }  j \in E:= \{0, \ldots, n \}^d \biggr\}
\]
can be described by the $d^{\text{th}}$ level RLT over $\Lambda$, which is denoted as $\RLT_d(\Lambda)$.
\begin{theorem}~\label{them:hull-Lambda}
$\conv(M^\Lambda) = \proj_{(\lambda,w)}\bigl(\RLT_d(\Lambda)\bigr) = R$, where 
\begin{equation*}
\begin{aligned}
R:=\biggl\{ (\lambda, w) \biggm|	 \;w \geq 0,\	\sum_{j \in E} w_j = 1,\  \lambda_{ij_i} = \sum_{p \in E:p_i =j_i}w_p \; \text{ for }  i \in [ d ] \text{ and }   j \in E \biggr\}. 
\end{aligned}
\end{equation*}
\end{theorem}
\begin{compositeproof}
First, we show  $\conv( M^\Lambda) = R$. Consider a set $M^{\vertex(\Lambda)}: = \bigl\{(\lambda, w)\in M^{\Lambda} \bigm| \lambda \in \vertex(\Lambda)\bigr\}$. Then, we obtain $ \conv(M^\Lambda) = \conv\bigl(M^{\vertex(\Lambda)}\bigr) = R$, where the first equality follows from Corollary 2.7 in~\cite{tawarmalani2010inclusion} and the second equality can be established by  using disjunctive programming~\cite{balas1998disjunctive}. Now, the proof is complete if we show $\proj_{(\lambda, w)}\bigl(\RLT_d(\Lambda)\bigr) \subseteq R$ since the $d^{\text{th}}$ level RLT over $\Lambda$ yields a convex relaxation for $M^\Lambda$. 
%For $I \subseteq \{1, \ldots, d\}$ and $p \in \{0, \ldots, n\}^d$, the monomial $ \prod_{i \in I}\lambda_{ip_i}$ is linearized by $y_{(I,p)}$. Aggregating $(1- \lambda_{i0} - \cdots - \lambda_{in})\prod_{i' \in I}\lambda_{ij_i}\prod_{i'\not\in I\cup\{i\}}(1-\sum_{j=0}^n \lambda_{i'j}) = 0$ over all $I\subseteq \{1,\ldots,d\}\backslash\{i\}$  and $(j_i')_{i'\ne i} \in \{0,\ldots,n\}^{d-1}$, we have $1-\sum_{j=0}^n \lambda_{ij} = 0$. Similarly, we can argue that $\lambda_{ij} \geq 0$ is implied at $d^{\text{th}}$ level RLT of $\Lambda$. 
For $I \subseteq [ d ]$ and $j \in E$, define $E(I,j) = \{p \in E \mid p_i = j_i \; \text{for } i \in I \}$. Then, the $d^{\text{th}}$ level RLT constraints, with $ \prod_{i \in I}\lambda_{ij_i}$ linearized by $y_{(I,j)}$, are
\[
y_{([ d ],j)}  \geq 0 \quad \text{and} \quad \sum_{I' \supseteq I, j' \in E(I,j)} (-1)^{|I'\setminus I|}y_{(I',j')} =0 \quad \text{for } I \subset [ d ] \text{ and } j \in E.
\]
 For $I \subseteq [ d ] $ and $j \in E$, let $w_{(I,j)}$ denote $ \prod_{i \in I}\lambda_{ij_i} \prod_{i \notin I}(1- \sum_{j=0}^n\lambda_{ij})$. Since, $y(I,j)$ and $w(I,j)$ each form a basis of the space of multilinear functions, it follows that there is an invertible transformation relating the two set of variables, that is 
 \[
 y_{(I,j)} = \sum_{I' \supseteq I, j' \in E(I,j)}w_{(I',j')} \quad \text{and} \quad w_{(I,j)} = \sum_{I' \supseteq I, j' \in E(I,j)} (-1)^{|I'\setminus I|}y_{(I',j')}.
 \]
Then, expressing $y_{(I,j)}$ in terms of $w_{(I,j)}$ variables, the $d^{\text{th}}$ level RLT relaxation can be represented by  constraints $w_{([ d ],j)} \geq 0$, $w_{(I,j)} = 0$ for $I \subset [ d ]$ and $j \in E$, and $y_{(I,j)} = \sum_{I' \supseteq I, j' \in E(I,j)}w_{(I',j')}$ for $I \subseteq [ d ]$  and $j \in E$. In particular, using the relation $ w_{([ d ],j)} = w_j$, $y_{(\{i\},j)} = \lambda_{ij_i}$ and $y_{(\emptyset,j)} = 1$, this system implies constraints in $R$. \Halmos
\end{compositeproof}

\revision{
Now, we apply Theorem~\ref{them:hull-Lambda}  to our setting. Consider a multilinear outer-function $\phi:\R^d\to \R$, that is $\phi(s_{1n}, \ldots, s_{dn}) = \sum_{I \in \mathcal{I}}c_I\prod_{i \in I}s_{in}$. Then, an explicit convex hull description of the graph of $\bar{\phi}(\cdot)$ over $Q$ is given as follows
\begin{equation}\label{eq:hull-Q-single}
\begin{aligned}
 \phi &= \sum_{I \in \mathcal{I}}  c_I\sum_{p \in E}\biggl(\prod_{i \in I} a_{ip_i}\biggr) w_p \quad w \geq 0 \quad \sum_{p \in E} w_p = 1,\\
 s_i &= \sum_{j = 0}^{n} v_{ij} \sum_{p \in E:p_i = j}w_p \quad \text{  for } i \in [d].
 \end{aligned}
\end{equation}
This holds since convexification commutes with affine transformation and the graph of $\bar{\phi}(\cdot)$ over $Q$ is the image of $M^\Lambda$ under an affine transformation given as follows:  $s_i = \sum_{j = 0}^nv_{ij} \lambda_{ij}$ for $i \in [ d ]$, where $\{v_{10}, \ldots, v_{in}\}$ is the vertex set of $Q_i$ defined in~(\ref{eq:Q-V}), $\phi = \sum_{I \in \mathcal{I}}c_Im_I$ and $m_{I}= \prod_{i \in I}s_{in} = \prod_{i\in I}s_{in} \bigl(\prod_{i\notin I}\sum_{p =0}^n \lambda_{ip}\bigr)= \sum_{p \in E} \bigl( \prod_{i \in I} a_{i  p_i }\bigr) w_{p}$ for each $I \in \mathcal{I}$. 	The convex hull description in~\eqref{eq:hull-Q-single}, together with Theorems~\ref{them:MIP-phi-z} and~\ref{them:DCR},  yields an explicit MIP relaxation for a composite function $\phi \mcirc f$ with a multilinear outer-function.
\begin{corollary}\label{cor:mul}
Consider a composite function $\phi \mcirc f$, where $\phi: \R^d \to \R$ is a multilinear function, that is $\phi(s_{1n}, \ldots, s_{dn}) = \sum_{I \in \mathcal{I}}c_I\prod_{i \in I}s_{in}$. An ideal formulation for $\{\conv(\graph(\phi|_H))\}_{H \in \mathcal{H}}$ is given by $\bigl\{(\phi, w,s,\delta) \bigm| (\phi, w, s) \in(\ref{eq:hull-Q-single}),\  Z(s) =z,\ (z, \delta) \in (\ref{eq:Inc-1}) \bigr\}$ which, together with  $u(x) \leq s$ and $(x,s_{\cdot n}) \in W$, yields an MICP relaxation for the graph of $\theta \mcirc f$. 
\end{corollary}
}

\section{Geometric insights and strengthening the relaxation}\label{section:comparison}
We first compare the strength of relaxations given by the composite relaxation scheme in~\cite{he2021new}, the disjunctive constraint in~(\ref{eq:MIP-standard}), and Theorem~\ref{them:DCR}.  Since they are not defined in the same space, we consider three functions $\varphi_{\mathcal{H}}(\cdot)$, $\varphi_{\mathcal{H}-}(\cdot)$, and $\varphi(\cdot)$ defined as follows:
\begin{equation*}\label{eq:comparison}
\begin{aligned}
	\varphi(x,s_{\cdot n}) &: = \max \Bigl\{ \conc_Q(\ephi)(s) \Bigm|    u(x) \leq s, (x,s_{\cdot n}) \in W \Bigr\},\\
		\varphi_{\mathcal{H}-}(x,s_{\cdot n}) &: = \max\Bigl\{ \conc_Q(\ephi)(s) \Bigm|  \bigl(Z(s), \delta \bigr) \in(\ref{eq:Inc-1}), (x,s_{\cdot n}) \in W \Bigr\},\\
	\varphi_{\mathcal{H}}(x,s_{\cdot n}) & := \max \Bigl\{ \conc_Q(\ephi)(s) \Bigm|
\bigl(Z(s),\delta \bigr) \in (\ref{eq:Inc-1}) ,  u(x) \leq s , (x,s_{\cdot n}) \in W
\Bigr\}.
\end{aligned}
\end{equation*}
The first function, $\varphi(\cdot)$, is derived using the composite relaxation defined as in~(\ref{eq:CR}). The second function, $\varphi_{\mathcal{H}-}(\cdot)$, uses Theorem~\ref{them:MIP-phi-z} to reformulate \eqref{eq:MIP-z} and relaxes the graph of inner functions as $W$. Similar relaxation ideas has appeared in the literature~\cite{misener2011apogee,misener2012global,nagarajan2019adaptive,huchette2019combinatorial}. However, as discussed in Example~\ref{ex:size-comparision}, these relaxations were either not ideal or require many continuous variables. In contrast, our definition of $\varphi_{\mathcal{H}-}(\cdot)$ exploits convexification results for the outer-function over $Q$ to obtain an ideal formulation that uses $d(n+1)$ continuous variables. The last function, $\varphi_{\mathcal{H}}(\cdot)$, is derived from Theorem~\ref{them:DCR}. Clearly, for every $(x,s_{\cdot n}) \in W$, $\varphi_{\mathcal{H}}(x,s_{\cdot n}) \le \min\bigl\{\varphi(x,s_{\cdot n}),\varphi_{\mathcal{H}-}(x,s_{\cdot, n})\bigr\}$, since in $\varphi_{\mathcal{H}}$, we relate an MIP formulation of $\{ \Delta_H\}_{H \in \mathcal{H}}$ with the underestimators of inner-functions via the inequalities $u(x)\le s$.
%As a result,  for every $(x,s_{\cdot n}) \in W$, $\varphi_{\mathcal{H}}(x,s_{\cdot n}) \leq \varphi_{\mathcal{H}-}(x, s_{\cdot n})$ and $\varphi_{\mathcal{H}}(x,s_{\cdot n}) \leq \varphi(x,s_{\cdot n})$. 
In the following, we derive geometric insights into relaxation quality, and use these insights to tighten the third function further in Theorem~\ref{them:DCR+}.

%More specifically, $\varphi_{\mathcal{H}-}(x,s_{\cdot n})$ is derived by using the reformulation \eqref{eq:MIP-z} of \eqref{eq:MIP-standard} and outer-approximating the inner functions using $(x,s_{\cdot n})\in W$. On the other hand, $\varphi_{\mathcal{H}}(x,s_{\cdot n})$ tightens this relaxation by, additionally appending the inequalities $u(x)\le s$ implicit in our construction~\eqref{eq:ordered-oa}.
%In other words, the MIP relaxation in Theorem~\ref{them:DCR} is tighter than the MIP relaxation given by the disjunctive constraints~(\ref{eq:MIP-standard}) and the composite relaxation in Proposition~\ref{prop:CR}. 
In Proposition~\ref{prop:eval}, we will show that the three formulations evaluate a non-increasing function at three different points. To derive this geometric insight, we introduce the following discrete function $\xi_{i,a_i}$ of $u_i$:
\begin{equation*}\label{eq:2-d-rep}
	\xi_{i,a_i}(a;u_i) = \left\{ \begin{aligned}
		&u_{ij} && a=a_{ij} \quad \text{for } j \in \{0\} \cup [ n ] \\
		&-\infty && \text{otherwise}. 
	\end{aligned}	
	\right.
\end{equation*}
For a given $\bar{x}$, assume that $u_{ij}(\bar{x})$ underestimates $\min\{a_{ij},f_i(x)\}$ and let
\[
s^{u}_i := \Bigl(\conc(\xi_{i,a_i})\bigl(a_{i0};u_i(\bar{x})\bigr), \ldots, \conc(\xi_{i,a_i})\bigl(a_{in};u_i(\bar{x})\bigl) \Bigr) \quad \text{for } i \in [ d ],
\]
 where each argument denotes the concave envelope of $\xi_{i,a_i}\bigl(\cdot;u_i(\bar{x})\bigr)$ over $[a_{i0}, a_{in}]$. Since $u_{ij}(\bar{x})$ underestimates $\min\{a_{ij},f_i(x)\}$, the constructed $s^{u}_i$ belongs to $Q_i$ (see Proposition 4 in~\cite{he2021new}). We show next that the three relaxations $\varphi(\cdot)$, $\varphi_{\mathcal{H}-}(\cdot)$, and $\varphi_{\mathcal{H}}(\cdot)$ use different $u(\cdot)$ to construct $s^{u}$, with that used by $\varphi_{\mathcal{H}}(\cdot)$ being the largest. Then, each relaxation evaluates $\conc_Q(\bar{\phi})(\cdot)$ at the corresponding $s^u$. Since $\xi_{i,a_i}$ is non-decreasing with $u_i(\bar{x})$, it follows that $s^{u}_i$ is non-decreasing with $u_i(\bar{x})$. Moreover, since $\conc_Q(\bar{\phi})$ is  a non-increasing function \cite{he2021new}, it follows that $\varphi_{\mathcal{H}}(\cdot)$ yields the tightest relaxation. The larger $u(\cdot)$ used by $\varphi_{\mathcal{H}}(\cdot)$ combines information from underestimators of the inner functions and the relevant discretization bounds to achieve its tightness.
\begin{proposition}\label{prop:eval}
Assume the same setup as in Theorem~\ref{them:DCR}. Let $(\bar{x},\bar{f}) \in W$, and let $\bar{t}:=(\bar{t}_1, \ldots, \bar{t}_d)$ is a vector such that $ \bar{f} \in \bar{H}:=\prod_{i=1}^d[a_{i \tau(i,\bar{t}_i-1) }, a_{i \tau(i,\bar{t}_i) } ]$. For $i \in [ d  ] $, consider vectors $\u_i$, $\hat{u}_i$ and $u^*_i$ in $\R^{n+1}$ defined as 
\begin{equation*}~\label{eq:updating-underestimation}
\u_{ij} = \begin{cases}
\bar{f}_i & j = n \\
	u_{ij}(\bar{x}) & j <n 
\end{cases},\ \hat{u}_{ij} = \begin{cases}
	a_{ij} & j \leq \tau(i,\bar{t}_i-1) \\
	\bar{f}_i& j \geq \tau(i,\bar{t}_i) \\
  -\infty & \text{otherwise}
\end{cases},\ \text{and } u^*_i =  \u_i \vee \hat{u}_i.
\end{equation*}
Then, we have  $\varphi(\bar{x},\bar{f})$, $\varphi_{\mathcal{H}-}(\bar{x},\bar{f})$, and $\varphi_{\mathcal{H}}(\bar{x},\bar{f})$ are $\conc_Q(\ephi)(\s)$, $\conc_Q(\ephi)(\hat{s})$, and $\conc_Q(\ephi)(s^*)$, respectively,
%\[
%\begin{aligned}
%	\varphi(\bar{x},\s_{\cdot n}) &= \conc_Q(\ephi)(\s) \\
%	\varphi_{\mathcal{H}-}(\bar{x},\s_{\cdot n}) &= \conc_Q(\ephi)(\hat{s}) \\
%	\varphi_{\mathcal{H}}(\bar{x},\s_{\cdot n}) &= \conc_Q(\ephi)(s^*)
%\end{aligned}
%\]
where for all $i$ and $j$, $\s_{ij} = \conc(\xi_{i,a_i})(a_{ij};\u_i)$, $\hat{s}_{ij}= \conc(\xi_{i,a_i})(a_{ij};\hat{u}_i)$ and $s^*_{ij}= \conc(\xi_{i,a_i})(a_{ij};u_i^*)$.
% Moreover, $\varphi(\bar{x},\bar{f}) \leq \varphi_{\mathcal{H}}(\bar{x},\bar{f})$ and $\varphi_{\mathcal{H}-}(\bar{x},\bar{f}) \leq \varphi_{\mathcal{H}}(\bar{x},\bar{f})$. 
\end{proposition}
\begin{compositeproof}
	See Appendix~\ref{app:eval}.
\end{compositeproof}

\newcommand{\BoundFunc}[1]{\setsepchar{:}\readlist*\ZZ{#1}{\cal L}_{\ZZ[1]\ZZ[2]}(\ZZ[3])}
\newcommand{\bb}[1]{\setsepchar{:}\readlist*\ZZ{#1}b_{\ZZ[1]\ZZ[2]\ZZ[3]}}
\newcommand{\bd}[1]{\setsepchar{:}\readlist*\ZZ{#1}\mathop{\text{bd}}_{\ZZ[1]\ZZ[2]\ZZ[3]}}
\newcommand{\idx}{k}
\newcommand{\jtoidx}[1]{\setsepchar{:}\readlist*\ZZ{#1}\theta(\ZZ[1],\ZZ[2])}
The remainder of this section derives Theorem~\ref{them:DCR+} which tightens the relaxation of Theorem~\ref{them:DCR}. Recall that
Theorem~\ref{them:DCR} yields a valid relaxation as long as $(u,a,f)$ satisfy \eqref{eq:ordered-oa}. To tighten the relaxation, we will use, for each underestimator, the best local bounds available when binary variables are fixed. Since these bounds are typically tighter than the global bounds, $a$, the resulting relaxation is also tighter than that of Theorem~\ref{them:DCR}. Towards this end, we define local bounds as a function of binary variables
\[
a'_{ij} = \BoundFunc{i:j:\delta_i} = a_{i0} + \sum_{\idx=0}^{l_i} \bb{i:j:\idx}\delta_{i\idx} \quad \text{ for } i \in [ d ]  \text{ and } j \in [ n ],
\]
where $b_{ijk}$ is an incremental bound, $\delta_{i0}$ is assumed to be one and $\delta_{il_i} = 0$. For notational convenience, we set $\tau(i,0)=0$. Here, $a_{i0}+b_{ij0}$ is the upper bound for $u_{ij}(x)$ over  $f_i(x)\in [a_{i0},a_{i\tau(i,1)}]$. For $1\le k< l_i$, $b_{ijk}$ is visualized as a difference of bounds. Let $\bar{a}_{ij}$ (resp. $\hat{a}_{ij}$) be the upper bounds for $u_{ij}(x)$ such that $f_i(x)\in [a_{i\tau(i,k-1)},a_{i\tau(i,k)}]$ (resp.  $f_i(x)\in [a_{i\tau(i,k)},a_{i\tau(i,k+1)}]$) then the incremental bound $b_{ijk}$ is $\hat{a}_{ij} - \bar{a}_{ij}$.
%By this, we mean that upper bounds for $u_{ij}(\cdot)$ are obtained iteratively. Assume that our upper bound for $u_{ij}(\cdot)$ over the region where $f_i(x)\in [a_{i\tau(i,k-2)},a_{i\tau(i,k-1)}]$ is $\bar{a}_{ij}$. Then, to obtain an upper bound over the region where $f_i(x)\in [a_{i\tau(i,k-2)},a_{i\tau(i,k-1)}]$, we add $b_{ijk}$ to $\bar{a}_{ij}$.
%That is, if the region in consideration is such that $f_i(x)$ exceeds $a_{i,\tau(i,k)}$ then $\delta_{ik}$ equals one. Then, to obtain a valid upper bound we add $b_{ijk}$ to the upper bound on $u_{ij}(\cdot)$ which was derived over the region $f_i(x)\le a_{i,\tau(i,k)}$.
%interpreted as the incremental bound that when added to the bound   $u_{ij}(x)$ as the binary variable $\delta_{ik}$ is activated.
In particular, for the setting of Theorem~\ref{them:DCR} we would define $\bb{i:j:0}=a_{ij}-a_{i0}$ and $\bb{i:j:k}=0$ for $k\ge 1$ because bounds were not updated based on the value of $f_i(x)$.  Our construction relies on two conditions on $\BoundFunc{i:j:\delta_i}$. First,  for $0\le j' \le j\le n$, we assume that for all $\idx \in \{0\} \cup [ l_i] $, $\sum_{\idx'=0}^\idx \bb{i:j:\idx'} \ge  \sum_{\idx'=0}^\idx \bb{i:j':\idx'}$. This condition imposes that for $\delta_i$ satisfying $1=\delta_{i0}\ge\cdots \ge\delta_{il_i}\ge \delta_{il_i+1}=0$, we have 
\begin{equation*}
    \BoundFunc{i:j:\delta_i} - \BoundFunc{i:j':\delta_i} = \sum_{\idx'=0}^{l_i}  \bigl(\delta_{i\idx'}-\delta_{i\idx'+1}\bigr)\biggl(\sum_{\idx''=0}^{\idx'} \bigl(\bb{i:j:\idx''}-\bb{i:j':\idx''}\bigr)\biggr) \ge 0.
\end{equation*}
It follows that $a'_{ij}\ge a'_{ij-1}$ for all $i \in [ d ] $ and $j\in [ n ]$. Second, we update local bounds in a particular way. We choose $b_{ijk}$ so that whenever $1=\delta_{ik-1}>\delta_{ik}=0$, {\it i.e.},\/ $f_i(x)\in [a_{i\tau(i,k-1)},a_{i\tau(i,k)}]$, $\BoundFunc{i:j:\delta_i}$ is between $a_{i,\tau(i,k-1)}$ and $a_{i,\tau(i,k)}$. More specifically,
\begin{alignat*}{2}
	&\BoundFunc{i:j:\delta_i} \geq \max\bigl\{ a_{i\tau(i,k-1)},\max_x\{u_{ij}(x)\mid f_i(x)\in [a_{i\tau(i,k-1)},a_{i\tau(i,k)}]\}\bigr\} &\quad&\text{ for } j < \tau(i,k)\\
	&\BoundFunc{i:j:\delta_i} = a_{i\tau(i,k)}  &&\text{ otherwise}.
\end{alignat*}
%Recall that $u_{ij}(\cdot)$ is an underestimator of $f(\cdot)$ bounded from above by $a_{ij}$. 

Recall that Theorem~\ref{them:DCR} made extensive use of the transformation~\eqref{eq:Z_inv}. Since the bounds depend on $\delta_i$, at first blush this transformation is bilinear. However, it simplifies to a linear expression for $s_i$ in terms of $(z_i,\delta_i)$ as follows. Letting $\jtoidx{i:j} = \arg\min\{\idx\mid \tau(i,\idx)\ge j\}$, if $(z_i,\delta_i)$ satisfies~(\ref{eq:Inc-1}), we rewrite \eqref{eq:Z_inv} as:
\begin{equation*}
\begin{split}
    s_{ij}-s_{ij-1} &=  z_{ij} \bigl(\BoundFunc{i:j:\delta_i}-\BoundFunc{i:j-1:\delta_i}\bigr)
    =\sum_{\idx=0}^{l_i} (\bb{i:j:\idx}-\bb{i:,j-1:,\idx})\delta_{i\idx}z_{ij} \\
    &= 
    z_{ij}\sum_{\idx=0}^{\jtoidx{i:j}-1} (\bb{i:j:\idx}-\bb{i:,j-1:,\idx}) + \sum_{\idx=\jtoidx{i:j}}^{l_i} (\bb{i:j:\idx}-\bb{i:,j-1:,\idx})\delta_{i\idx},
\end{split}
\end{equation*}
where the simplification in the last equality occurs since $\delta_{i\idx}z_{ij} = z_{ij}$ if $\idx < \jtoidx{i:j}$ and $\delta_{i\idx}$ otherwise. We define $s_{i0} = a'_{i0} = a_{i0} + \sum_{\idx=0}^{l_i} b_{i0\idx}\delta_{i\idx}$ so that
\begin{equation*}
  \begin{split}
  s_{ij} &= a_{i0} + \sum_{\idx=0}^{l_i} \bb{i:0:\idx}\delta_{i\idx} + \sum_{j'=1}^j \biggl\{ z_{ij'}\sum_{\idx=0}^{\jtoidx{i:j'}-1} (\bb{i:j':\idx}-\bb{i:,j'-1:,\idx}) \\
  &\qquad\qquad + \sum_{\idx=\jtoidx{i:j'}}^{l_i} (\bb{i:j':\idx}-\bb{i:,j'-1:,\idx})\delta_{i\idx}\biggr\}\\
  &=a_{i0} + \sum_{j'=1}^j z_{ij'}\sum_{\idx=0}^{\jtoidx{i:j'}-1} (\bb{i:j':\idx}-\bb{i:j'-1:\idx})
  +\sum_{\idx=0}^{l_i}\delta_{i\idx}\bb{i:,\min\{\tau(i,\idx),j\}:,\idx}\\
  &=a_{i0} + \sum_{j'=1}^j z_{ij'}\sum_{\idx=0}^{\jtoidx{i:j'}-1} (\bb{i:j':\idx}-\bb{i:j'-1:\idx})\\
  &\qquad\qquad+\sum_{\idx=0}^{\jtoidx{i:j}-1}\delta_{i\idx}\bb{i:,\tau(i,\idx):,\idx} +\sum_{\idx=\jtoidx{i:j}}^{l_i}\delta_{i\idx}\bb{i:,j:,\idx}:=G_{ij}(z,\delta)
  \end{split}
\end{equation*}
Observe that, in the above expression, the coefficients of $z_{ij'}$ for all $j'$ are non-negative because we have assumed that $\sum_{\idx=0}^{\jtoidx{i:j'}-1} (\bb{i:j':\idx}-\bb{i:j'-1:\idx})\ge 0$. If we additionally assume that $\bb{i:j:\idx}\ge 0$ for all $i\in [ d ]$, $j\in [ n ]$ and $ \idx \in \{0\} \cup [ l_i ] $ then all the coefficients in the above definition of $s_{ij}$ are non-negative. Let $\delta$ be binary so that for each $i$, there is a $\idx_i$ so that $1 = \delta_{i\idx_i-1} > \delta_{i\idx_i}=0$. Then, it follows that $a'_{ij}=\BoundFunc{i:j:\delta_i} = a_{i0} +\sum_{\idx=0}^{{\idx_i}-1}\bb{i:j:\idx}$ and the above calculation shows that $s_{ij} = a'_{i0} + \sum_{j'=1}^j z_{ij'}(a'_{ij'}-a'_{ij'-1})$. In other words, $s_i$ belongs the simplex $Q'_i$ defined as in~(\ref{eq:Q-V}) using vector $a'_i = (a'_{i0},\ldots, a'_{in})$.

\begin{theorem}~\label{them:DCR+}
 Assume the same setup as in Theorem~\ref{them:DCR}. Then, an MICP relaxation for the hypograph of $\phi \mcirc f$ is given by
\begin{equation}\label{eq:DCR+}
\left\{ (x, \phi, s, z, \delta) \left| \;
\begin{aligned}
&\phi \leq \conc_{\hat{\Delta}}(\ephi \mcirc G )(z,\delta),\ s= G(z,\delta),\ (z,\delta ) \in (\ref{eq:Inc-1}) \\
&u(x) \leq s,\ (x,s_{\cdot n}) \in W	
\end{aligned}
\right.
\right\},
\end{equation}
where $\hat{\Delta}$ is the simplotope in $(z,\delta)$ variables given by~(\ref{eq:Inc-1}).  
\end{theorem}
\begin{compositeproof}
See Appendix~\ref{app:DCR+}.
\end{compositeproof}
\begin{remark}
We remark that the relaxation in Theorem~\ref{them:DCR+} improves that the one given by Theorem~\ref{them:DCR}. Consider a function $\phi_{\mathcal{H}_+}(x,s_{\cdot n}): = \max\bigl\{ \phi \bigm| (x,\phi,s,z,\delta) \in(\ref{eq:DCR+}) \bigr\}$. It can be shown that for $(\bar{x},\bar{f}) \in W$, $\phi_{{\cal H}}(\bar{x},\bar{f}) \geq \phi_{{\cal H}_+}(\bar{x},\bar{f}) =\conc_{Q'}(\bar{\phi})(s')$,  where $Q'$ is formed using $a'$, and $s'_{ij} = \conc(\xi_{i,a'_i})(a'_{ij};\tilde{u}_i)$, where $\tilde{u}_i$ is the same as $u^*_i$ in Proposition~\ref{prop:eval} except that $a_{i}$ is replaced with $a'_{i}$ (See Appendix 11 of~\cite{he2021new}).  \QED
\end{remark}
                                                                                                                                                                                                                    
\begin{remark}
Theorem~\ref{them:DCR+} requires the concave envelope of $\ephi\mcirc G$ over $\hat{\Delta}$. In this remark, we show that this envelope is readily available under certain conditions. Assume $\phi(\cdot)$ is supermodular and $\bar{\phi}\circ G$ is concave extendable from $\vertex(\hat{\Delta})$. The latter condition would be satisfied for a multilinear $\phi$. We provide an explicit description of the concave envelope  $\ephi \mcirc G$ over $\hat{\Delta}$ for this case.  First, observe that $\vertex(\hat{\Delta}) = \prod_{i=1}^d \vertex(\hat{\Delta}_i)$, and $\vertex(\hat{\Delta}_i)$ forms a chain with join (resp. meet) defined as component-wise maximum (resp. minimum). Moreover, if $\phi(\cdot)$ is supermodular over $[f^L,f^U]$ then $\ephi \mcirc G$ is supermodular when restricted to $\vertex(\hat{\Delta})$, \textit{i.e.}, for $y':=(z',\delta')$ and $y'':=(z'',\delta'')$ in $\vertex(\hat{\Delta})$, 
\[
\begin{aligned}
(\ephi \mcirc G)(y' \vee y'')  &+(\ephi \mcirc G)(y' \wedge y'')  =  \phi \bigl( G_{1n}(y'_1) \vee G_{1n}(y''_1), \ldots,G_{dn}(y'_d) \vee G_{dn}(y''_d) \bigr) + \\
& \qquad \qquad \qquad \phi \bigl( G_{1n}(y'_1) \wedge G_{1n}(y''_1), \ldots,G_{dn}(y'_d) \wedge G_{dn}(y''_d) \bigr) \\
&\geq \phi \bigl(G_{\cdot n}(y') \bigr) + \phi \bigl(G_{\cdot n}(y'') \bigr) = (\ephi \mcirc G)(y') + (\ephi \mcirc G)(y''),
\end{aligned}
\]
where the first equality holds since $G_{in}$ is non-decreasing and $\vertex(\hat{\Delta}_i)$ forms a chain, the inequality holds by supermodularity of $\phi(\cdot)$, and the last equality holds by definition. Therefore, if $ \ephi \mcirc G$ is concave-extendable from $\vertex(\hat{\Delta})$ and $\phi(\cdot)$ is supermodular on $[f^L,f^U]$, the concave envelope of $\ephi \mcirc G$ over $\hat{\Delta}$ is given by Proposition 3 in~\cite{he2022tractable}.
%-- a consequence of Corollary 3.4 in~\cite{tawarmalani2013explicit}.
\QED
\end{remark}

\section{Computational experiments}\label{section:computation}
This section describes our preliminary implementation of discrete and continuous relaxations for MINLP, and provides numerical comparisons on their tightness. 

\subsection{Implementation via a propagation algorithm}\label{section:algorithm}
\renewcommand{\algorithmicrequire}{\textbf{Input:}}
\renewcommand{\algorithmicensure}{\textbf{Output:}}

\begin{algorithm}
\caption{Relaxation Propagation Procedure} \label{alg:propagation}
	\begin{algorithmic}[1]
	\Require An expression tree representation $T$ of a factorable function.
	\Ensure An array of tuples $\textsc{Rlx}(T)=(u(\cdot),a,\chi)$ for root node that satisfies~(\ref{eq:ordered-oa-over})
		\If{the root node expression $v(\cdot)$ of $T$ is multiplication}
			\State $T_1 \gets$ the left subtree
			\State $T_2 \gets$ the right subtree
			\For{$(u_{1}(\cdot),a_{1}, \chi_1)  \in \textsc{Rlx}(T_1)$ and $(u_{2}(\cdot),a_{2}, \chi_2)  \in \textsc{Rlx}(T_2) $}
			\State introduce variables $(s_1, f_1, s_2, f_2)$
			\For{$i = 1,2$}
						\If{$\chi_i$ is true}
			\State add constraints $u(\cdot)_i \leq s_i \leq a_i$ and $s_i \leq f_i$
			\Else 
			\State add constraints $a_i \leq s_i \leq u(\cdot)_i$ and $f_i \leq s_i$
			\EndIf
			\EndFor
			\State \revision{apply~\eqref{eq:12-inequalities} to generate  $e_j(\cdot)\leq v(\cdot)$ and $v(\cdot)\leq r_k(\cdot)$, where $j,k \in [ 6 ]$}
				\For{ $j = 1, \ldots,6$}
				\State find $a$ such that $e_i(\cdot) \leq a$, and push $(e_i(\cdot), a, \text{true})$ to $\textsc{Rlx}(T)$
				\State find $a$ such that $a \leq r_k(\cdot)$, and push $(r_k(\cdot), a, \text{false})$ to $\textsc{Rlx}(T)$
				\EndFor
			\EndFor
			\ElsIf{the root node expression $v(\cdot)$ is a unary function}  
			\State $\hat{v}(\cdot)\gets$ the concave envelope of $v$ over bounds $[f^L,f^U]$ on its child node $f$
			\State $\check{v}(\cdot)\gets$ the convex envelope of $v$ over bounds $[f^L,f^U]$ on its child node $f$
			\For{a list of points $p^j$ in $[f^L,f^U]$}
			\State find the subgradient inequality $e_j(\cdot) \leq \check{v}(\cdot)$ at point $p^j$ and $a$ such that $e_i(\cdot) \leq a$, and push $(e_i(\cdot), a, \text{true})$ to $\textsc{Rlx}(T)$
			\State find the supergradient inequality $\hat{v}(\cdot) \leq r_j(\cdot)$ at point $p^j$ and $a$ such that $a \leq r_k(\cdot)$, and push $(r_k(\cdot), a, \text{false})$ to $\textsc{Rlx}(T)$
			\EndFor
			\EndIf
			\State Sort $\textsc{Rlx}(T)$ by the second coordinate ($a$) of each tuple so that~(\ref{eq:ordered-oa-over}) is satisfied. 

%		\For{$i$ from $1$ to $d$}\label{alg:sign:for:loop:dimensions}
%		\For{$j$ from $1$ to $n-1$}
%		\If{$\alpha_{ij} > 0$}
%		\State $\text{push}(j, J^+_i)$\label{alg:sign:initialize:pos}
%		\EndIf
%				\EndFor
%				\For{$j$ from $0$ to $n$}
%				\State $\prev(j) = j-1$\label{alg:sign:set:prev}
%				\State $\succ(j) = j+1$\label{alg:sign:set:succ}
%				\State $\alpha'_{ij} = \alpha_{ij}$;
%				\EndFor\label{alg:sign:initialize:done}
%			\While{$  J_i^+ \neq \emptyset $ }\label{alg:sign:while}
%			\State $j = \text{pop}(J_i^+)$\label{alg:sign:step:pickj}
%			\State $ = \frac{a_{i\succ(j)} - a_{ij}}{a_{i\succ(j)}-a_{i\prev(j)}}$
%				\For {$(j',\upmult)$ in $\bigl[\bigl(\prev(j),\bigr), \bigl(\succ(j),1-\bigr)\bigr]$}\label{alg:sign:step:for:loop:begins}
%			  	\State $\neg \gets (\alpha_{ij'} \le 0)$;
%					\State $\alpha'_{ij'} = \alpha'_{ij'} + \upmult\alpha'_{ij}$;\label{alg:sign:step:update_alpha}
%			  	\If{$\neg = \true$ and $\alpha'_{ij'} > 0$ and $0 < j' < n$}
%					\State $\text{push}(j',J_i^+)$\label{alg:sign:add:j'}
%				  \EndIf
%					\EndFor\label{alg:sign:step:for:loop:ends}
%				\State $\prev\bigl(\succ(j)\bigr) = \prev(j)$; $\succ\bigl(\prev(j)\bigr) = \succ(j)$;\label{alg:sign:remove:j}
%				\State $\prev(j) = -1$; $\succ(j) = n+1$;\label{alg:sign:orphan:j}
%				\State $\alpha'_{ij} = 0$;\label{alg:sign:zero:alpha}
%       \EndWhile
%     \EndFor
%		 \State \textbf{return} $\alpha'$.
		
	\end{algorithmic}
\end{algorithm}

For our computations, we focus on factorable programs, where each function is expressed as a recursive sum and product of univariate functions. For sums, we construct separate relaxations for each of the children. Consider now a node that represents the product of its child nodes, say $w=f_1(x)f_2(x)$. Iteratively, we select a pair $(u_1(\cdot),u_2(\cdot))$ of estimators for $(f_1(\cdot),f_2(\cdot))$, and obtain $(a_1,a_2)$ that satisfy $(u_1(x),u_2(x)) \le (a_1,a_2)$ over $X$. Assume that $f(x)$ is bounded so that $f^L \le f(x) \le f^U$. We introduce new variables $f_i$ and $s_i$, where variable $f_i$ represents $f_i(\cdot)$ and variable $s_i$ dominates $u_i(\cdot)$. \revision{Then, we relax the product $f_1f_2$ using the following $12$ inequalities:}  
%\begin{comment}
\begin{equation}~\label{eq:12-inequalities}
\begin{aligned}
f_1f_2 &\ge \max\left\{
			\begin{aligned}
				&e_1:=f_2^Uf_1 + f_1^Uf_2 - f_1^Uf_2^U \\
				&e_2:=(f^U_2-a_2)s_1 + (f^U_1-a_1)s_2 + a_2f_1 + a_1f_2  \\
				& \qquad \qquad+ a_1a_2- a_1f^U_2 - f^U_1a_2\\
				&e_3:=(f^U_2-f^L_2)s_1 + f^L_2f_1 + a_1f_2 - a_1f^U_2 \\
				&e_4:=(f^U_1-f^L_1)s_2 + a_2f_1 + f^L_1f_2 - f^U_1a_2 \\
				&e_5:= (a_2-f^L_2)s_1 + (a_1-f^L_1)s_2 + f^L_2f_1 + f^L_1f_2 - a_1a_2\\
				&e_6:=f_2^Lf_1 + f_1^Lf_2 - f_1^Lf_2^L \\
		\end{aligned}\right. \\
f_1f_2 &\le \min\left\{
			\begin{aligned}
				&r_1:= f_2^Uf_1 + f_1^Lf_2 - f_1^Lf_2^U \\
				&r_2:=  (f^L_2-a_2)s_1 + (a_1-f^U_1)s_2 + a_2f_1 + f^U_1f_2 - a_1f^L_2\\
				&r_3:=  (f^L_2-f^U_2)s_1 + a_1f_2 + f^U_2 f_1 - a_1f^L_2 \\
				&r_4:=  (f^L_1-f^U_1)s_2 + a_2f_1 + f^U_1 f_2 - f^L_1a_2 \\
				&r_5:=  (a_2-f^U_2)s_1 + (f^L_1-a_1)s_2 + f^U_2f_1 + a_1f_2 - f^L_1a_2\\
				&r_6:= f_2^Lf_1 + f_1^Uf_2 - f_1^Uf_2^L,
		\end{aligned}\right.
		\end{aligned}
		\end{equation}
%\end{comment}
where  $\max\{e_1,e_6\}\le f_1f_2\le \min\{r_1,r_6\}$ are McCormick inequalities (see \cite{mccormick1976computability}). These inequalities follow from Corollary~\ref{cor:bilinear-stair} and a direct proof of their validity is given in Theorems 1 and 5 of~\cite{he2021new}.
We also use the polynomial time procedure \revision{mentioned in} Remark~\ref{rmk:stair-Q} to separate inequalities, which we do not propagate to parent nodes. 

If $f_i(\cdot)$ has a non-trivial overestimator $o_i(\cdot)$ that is bounded from below by $a_i$, the inequalities are derived by replacing $o_i(x)$ with the underestimator $u_i(x) = f_i(x)-o_i(x)+a_i$, as discussed in Remark~\ref{rmk:over-and-under}.  Univariate functions are relaxed using subgradients (resp. supergradients) or their convex (resp. concave) envelopes. 
%Inductively, we construct relaxations for each of its children. These relaxations yield underestimators and overestimators for the function represented by each child node. 
The bounds for all linear estimators constructed are obtained  by simply maximizing and minimizing each term using its bounds. We next illustrate, on an example setting, the main ideas in propagation procedure of Algorithm~\ref{alg:propagation}.
\begin{example}\label{ex:recursive}
Consider $x_1^2x_2^2x_3^2$ over  $[1,2]^3$, and the expression tree in Figure~\ref{fig:expression-tree-1}, where edges are labeled with bounds of the tail nodes. Since $x_i^2\in [1,4]$ when $x_i\in [1,2]$, we obtain the convex underestimator using (factorable programming) FP, $x_1^2x_2^2 \geq \max \bigl\{x_1^2 + x_2^2 - 1, 4x_1^2 + 4x_2^2 - 16 \bigr\}$.
%\begin{equation*}\label{eq:recursive}
%x_1^2x_2^2 \geq \max \bigl\{x_1^2 + x_2^2 - 1, 4x_1^2 + 4x_2^2 - 16 \bigr\}.
%\end{equation*}
Then, inferring that $x_1^2x_2^2 \in [1,16]$, FP recursively underestimates the root node using $\max \{x_1^2 + x_2^2 - 1, 4x_1^2 + 4x_2^2 - 16 \} + x_3^2 -1$ and $4 \max \bigl\{x_1^2 + x_2^2 - 1, 4x_1^2 + 4x_2^2 - 16 \bigr\} + 16 x_3^2 -64$.
%\[
%x_1^2x_2^2x_3^2 \geq \max
%\left\{
%\begin{aligned}
%&\max \{x_1^2 + x_2^2 - 1, 4x_1^2 + 4x_2^2 - 16 \} + x_3^2 -1 \\ 
%&4 \max \bigl\{x_1^2 + x_2^2 - 1, 4x_1^2 + 4x_2^2 - 16 \bigr\} + 16 x_3^2 -64
%\end{aligned} 
%\right\}.
%\]
%Observe that $e(x_1,x_2) = 4x_1 +4x_2 - 9$ underestimates $x_1^2x_2^2$ and is upper bounded by $7$.
% Using Corollary~\ref{cor:bilinear-stair} we have $ x_1^2x_2^2x_3^2 \geq 3 e(x_1,x_2) + x_1^2x_2^2  +   7x_3^2 - 28 \geq 3 e(x_1,x_2) + 4x_1^2 + 4x_2^2 + 7x_3^2-44$, %Corollary~\ref{cor:stair-switching} 
%where the second inequality is obtained by further underestimating $x_1^2x_2^2$ using $4x_1^2+4x_2^2 -16$. At various points in $[1,2]^3$, this convex underestimator strictly dominates the factorable relaxation. 
%This shows that information about underestimators for lower-levels of the tree can be utilized while relaxing the higher-levels. 
%Now, exploiting the following ordering relation obtained from the left child of the root node
%\[
%x_1^2 + x_2^2 - 1 \leq \max \bigl\{x_1^2 + x_2^2 - 1, 4x_1^2 + 4x_2^2 - 16 \bigr\}  \quad \text{and} \quad x_1^2 + x_2^2 - 1 \leq 7 \quad \text{for } x \in [1,2]^3,
%\]
%we will construct a convex underestimator $3 e_1(x) +e_6(x)  +   7u_{32}(x) - 28$ $16x_1^2 + 16x_2^2 +  7x_3^2 - 64$ for $x_1^2x_2^2x_3^2$, which is not implied by the factorable one. 
%
Consider an expression tree depicted in Figure~\ref{fig:expression-tree-2}, whose edges are labeled with underestimators of tail nodes and their upper bounds over the box domain $[1,2]^3$. Namely, for $i \in \{1, 2, 3\}$, let $u_i(x):=(1,2x_i-1, 4x_i-4)$ and $a_i := (1,3,4)$. Exploiting the ordering relation, $u_{ij}(x) \leq x_i^2$, $u_{i2}(x) = x_i^2$, and $u_i(x) \leq a_i$ for $i \in \{1,2\}$, we obtain a relaxation for the left child of the root node:
\[
x_1^2x_2^2 \geq \max \left \{
\begin{alignedat}{3}
&e_0 := 1,  &\;&e_1:= s_{12} + s_{22}- 1\\
&e_2:=2s_{11}+ s_{12}  +  2s_{21} + s_{22} - 9,
&&e_3:=3s_{11} + s_{12}+ 3s_{22} - 12\\
&e_4:=3s_{12} + 3s_{21} + s_{22} - 12,
&&e_5:=s_{11} + 3s_{12} +s_{21} + 3s_{22} - 15\\
&\omit\rlap{$\displaystyle e_6:= \max \bigl\{e_0, \ldots, e_5, 4s_{12}+ 4s_{22} -16\bigr\}$}
\end{alignedat}
\right\}
,
\]
where $(1,2x_i-1,\max\{u_{ij}(x_i)\mid j=1,2,3\})\leq s_{i\cdot}\le (1,3,4)$ and then derive upper bounds for $e_i$ as $(1, 7,11,13,13,15,16)$. Then, these inequalities are used recursively to obtain estimators for the root node.
\QED
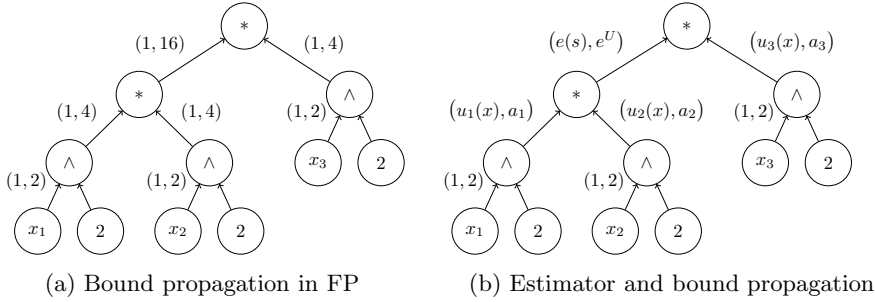
\begin{figure}[h]
\begin{subfigure}{0.48\linewidth}
	\centering
	\begin{tikzpicture}[scale=0.7, every tree node/.style={draw,circle,minimum size=2.5em}, level distance = 1.3cm, sibling distance=0.3cm,
edge from parent/.style={draw, <-,edge from parent path={(\tikzparentnode) -- (\tikzchildnode)}}]
\Tree[ .{$*$} \edge node[auto=right,pos=.5] {$(1,16)$};[.{ $ *$}  \edge node[auto=right,pos=.5] {$(1,4)$};[.$\wedge$  \edge node[auto=right,pos=.5] {$(1,2)$};$x_1$ 2 ]  \edge node[auto=left,pos=.5] {$(1,4)$};[.$\wedge$ \edge node[auto=right,pos=.5] {$(1,2)$};$x_2$ 2 ] ]  \edge node[auto=left,pos=.5] {$(1,4)$};[. $\wedge $  \edge node[auto=right,pos=.5] {$(1,2)$};$x_3$   $2$ ] ]
\end{tikzpicture}
\caption{Bound propagation in FP}\label{fig:expression-tree-1}
\end{subfigure}
\begin{subfigure}{0.48\linewidth}
 \begin{tikzpicture}[scale=0.7, every tree node/.style={draw,circle,minimum size=2.5em}, level distance = 1.3cm, sibling distance=0.3cm,
edge from parent/.style={draw, <-,edge from parent path={(\tikzparentnode) -- (\tikzchildnode)}}]
\Tree[ .{$*$} \edge node[auto=right,pos=.5] {$\bigl(e(s),e^U\bigr)$};[.{ $ *$}  \edge node[auto=right,pos=.5] {$\bigl(u_1(x),a_1\bigr)$};[.$\wedge$  \edge node[auto=right,pos=.5] {$(1,2)$};$x_1$ 2 ]  \edge node[auto=left,pos=.5] {$\bigl(u_2(x),a_2\bigr)$};[.$\wedge$ \edge node[auto=right,pos=.5] {$(1,2)$};$x_2$ 2 ] ]  \edge node[auto=left,pos=.5] {$\bigl(u_3(x),a_3\bigr)$};[. $\wedge $   \edge node[auto=right,pos=.5] {$(1,2)$};$x_3$   $2$ ] ]
\end{tikzpicture}
\caption{Estimator and bound propagation}\label{fig:expression-tree-2}
\end{subfigure}
\caption{Illustration of two propagation strategies on an expression tree.}
\end{figure}
%Observe that $u_{ij}(\cdot)$ can be substituted with their defining relations to obtain inequalities that do not require introduction of variables  beyond factorable programming scheme.\Halmos

%For notational simplicity, we considered the product of $x_1^2$, $x_2^2$ and $x_3^2$ but the construction generalizes naturally to $f(x)g(y)h(z)$, and as described before the example, to arbitrary expression trees, if each node can be transformed to be supermodular and/or submodular.
%\[
%x_1^2x_2^2x_3^3 \geq \max \left\{
%\begin{aligned}
%&	0 \\
%&	3e_2(x) + 9u_3(x) - 27 \\
%&	3 \max\bigl\{ e_3(x),e_4(x)\bigr\} + 12u_3(x) - 36 \\
%&	3e_5(x) + 15u_3(x) -45 \\
%&	3e_6(x) + 16u_3(x) - 48 \\
%&	4e_2(x) + 9u_3(x) - 36 \\
%&	4 \max\bigl\{ e_3(x),e_4(x)\bigr\} + 12u_3(x) - 48 \\
%&	4e_5(x) + 15u_3(x) -60 \\
%&	e_2(x)+3 \max\bigl\{ e_3(x), e_4(x) \bigr\}+3u_3(x)+9f_3(x)-45 \\
%&	e_2(x)+3e_5(x)+6u_3(x)+9f_3(x)-54\\
%&	e_3(x)+3e_5(x)+3u_3(x)+12f_3(x)-57\\
%&	e_2(x)+3e_6(x)+7u_3(x)+9f_3(x)-57\\
%&	e_3(x)+3e_6(x)+4u_3(x)+12f_3(x)-60\\
%&	e_5(x)+3e_6(x)+u_3(x)+15f_3(x)-63 \\
%&	4e_6(x) + 16f_3(x) - 64 
%\end{aligned}
%\right\}
%\]
\end{example}

\subsection{Numerical results}~\label{section:numerical}
Consider polynomial optimization problems:
\begin{equation}~\label{eq:rand-poly}
	\begin{aligned}
		\min\Bigl\{ \langle c, x \rangle + \langle d, y \rangle\Bigm|Ax + B y \leq b ,\ x^L \leq x \leq x^U ,\ y = (x^{\alpha_1}, x^{\alpha_2}, \ldots, x^{\alpha_m}) \Bigr\},
	\end{aligned}
\end{equation}
where $x\in \R^n$, $x^L\le x\le x^U$, $x^{\alpha_j} = \prod_{i=1}^n x_{i}^{\alpha_{ji}}$, $\alpha_j = (\alpha_{j1},\ldots,\alpha_{jn})$, $c\in \R^n$, $d\in\R^m_+$, $A\in \R^{r\times n}$, $B\in \R^{r\times m}_+$, and  $b\in \R^q$. We describe problem sizes by the parameter vector $(n,m,r)$. For each $(n,m,r) \in \{(15,30,10), (25,50,10), (50,100,20), (100,200,20)\}$,  we randomly generate $100$ instances of~(\ref{eq:rand-poly}) using a procedure detailed in Appendix~\ref{app:generation}. Our computational experiments are performed on a MacBook Pro with Apple M1 Pro with 10-cores CPU and 16 GB of memory. The code is written in \textsc{Julia v1.6} and our relaxations are modeled using \revision{\textsc{JuMP v1.4.0}~\cite{lubin2023jump}} and solved using \revision{\textsc{Gurobi v10.0.2}~\cite{gurobi}} as an MIP solver. We use \textsc{Ipopt v3.14}~\cite{wachter2006implementation} as an NLP solver to obtain upper bounds, and \revision{\textsc{SCIP v8.0}~\cite{bestuzheva2023global}} as an MINLP solver. 

We report on the quality of relaxations with five algorithmic settings.
\begin{itemize}
\item \texttt{SCIP}: The root node relaxation generated by \textsc{SCIP}, where we use the default settings, allowing cuts and range-reductions. \revision{The reported bound  is obtained just before the solver partitions the root node}. 
%We use \textsc{SCIP} for comparison because it is a state-of-the-art solver based on the factorable programming scheme. Our methods are designed to improve this scheme without introducing many new variables so that our relaxations scale to larger problem instances.
\end{itemize} 
For the next four settings, we use binary expression trees for monomials in (\ref{eq:rand-poly}) where leaf nodes are power functions and non-leaf nodes represent products of two subexpressions. The power functions are relaxed using sub-gradient inequalities at five uniformly selected points. For MIP relaxations,  we add a single discretization point at the median of the underestimator bounds. We do not match any common sub-expressions except for the univariate functions. The procedures below differ in how non-leaf nodes are relaxed.
\begin{itemize}
	\item[$\bullet$] \texttt{MC}: Each product is relaxed by replacing it with the four McCormick envelopes \cite{mccormick1976computability} following the FP scheme~\cite{tawarmalani2004global,belotti2009branching,misener2014antigone,vigerske2018scip,mahajan2021minotaur}.
    \item [$\bullet$] \texttt{MIP19}: \revision{Each product is relaxed by replacing it with the MIP relaxation obtained using formulation techniques in~\cite{huchette2019combinatorial}.}
	\item[$\bullet$] \texttt{MIP}: This relaxation is based on $\varphi_{\mathcal{H}^-}$ described in Section~\ref{section:comparison}. 
	\item[$\bullet$] \texttt{CR}: \revision{The composite relaxation is based on $\varphi$ described in Section~\ref{section:comparison}. Our implementation is based on the propagation procedure Algorithm~\ref{alg:propagation} and  Corollary~\ref{cor:bilinear-stair}. 	The inequalities in Corollary~\ref{cor:bilinear-stair} are implemented using a separation oracle, and we limit the number of separation calls to $50$.}
	\item[$\bullet$] \texttt{CRMIP:} \revision{This relaxation is based on $\varphi_{\mathcal{H}}$ described in Section~\ref{section:comparison}. Our implementation is   similar to that in \texttt{CR} except that, as in $\varphi_{\mathcal{H}}$, we augment \texttt{CR} with an incremental formulation~\eqref{eq:Inc}.} 
%	\item[$\bullet$] \texttt{CRMIP+}	
\end{itemize}
Relaxations similar to \texttt{MIP} have been proposed before \cite{misener2011apogee,misener2012global,nagarajan2019adaptive,huchette2019combinatorial}. However, \texttt{MIP} differs since it requires fewer continuous variables and uses limited separation rounds. \revision{Nonlinear JuMP models of all instances and sample \texttt{MIP} and \texttt{CRMIP} formulations are uploaded to zenodo~\cite{dataset}.}

%Since we do not implement a callback to separate the inequalities in Corollary~\ref{cor:bilinear-stair}, the corresponding relaxation can be weaker but scales more easily to large problem instances.

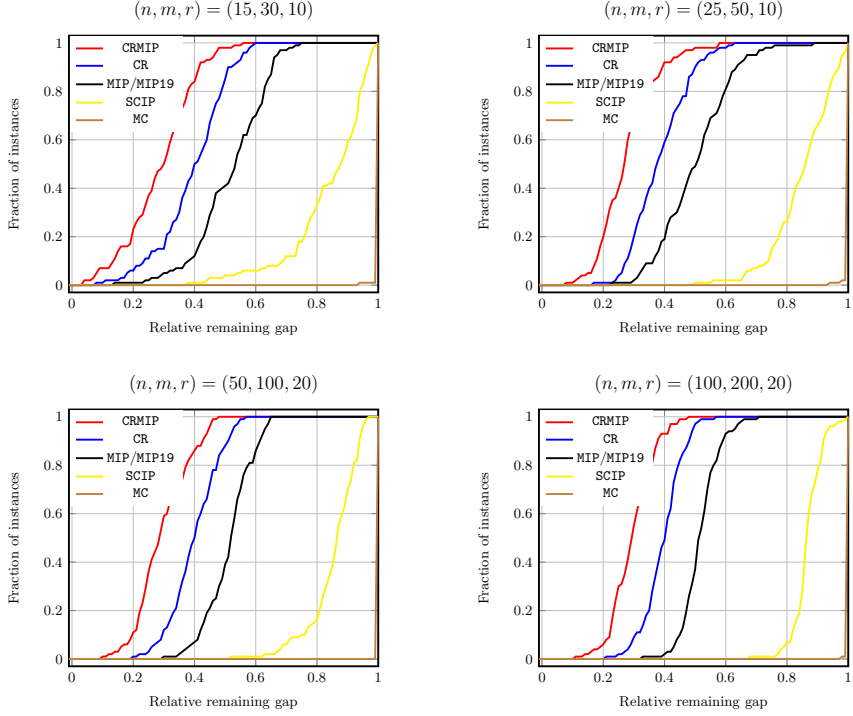
\begin{figure}[ht]
\begin{subfigure}{0.48\linewidth}
\begin{tikzpicture}[scale=0.6]
\pgfplotsset{
    tick label style={font=\small},
    label style={font=\small},
    legend style={font=\small},
    every axis/.append style={
    very thick,
    tick style={semithick}}
    }
\begin{axis}[
	title={\large $(n,m,r) = (15,30,10)$},
	xmin=-0.01,   xmax=1,
	ymin=-0.03,   ymax=1.03,
	xlabel={Relative remaining gap},
	ylabel={Fraction of instances},
	xtick={0,0.2,0.4,0.6,0.8,1},
	ytick={0,0.2,0.4,0.6,0.8,1}, 
	legend style={draw=none,
        at={(0.37,1)}},
    grid=major,
%   very thick,
    cycle list name=color list
        ]
\pgfplotstableread{size1.txt}\mydata;
\addplot table [x=poi,y=crmip] {\mydata};
\addplot table [x=poi,y=cr] {\mydata};
\addplot table [x=poi,y=mip] {\mydata};
%\addplot table [x=poi,y=mip1] {\mydata};
\addplot table [x=poi,y=scip] {\mydata};
\addplot table [x=poi,y=mc] {\mydata};
 \legend{$ \texttt{CRMIP}$, $\texttt{CR}$,$\texttt{MIP}/\texttt{MIP19}$,$\texttt{SCIP}$, $\texttt{MC}$}
%\addplot[smooth, point meta=explicit] coordinates {
%	(-2,3) (-1.5,2) (-0.3,-0.2) 
%	(1,1.2) (2,2) (3,5)};
\end{axis}
\end{tikzpicture}
\end{subfigure}
\begin{subfigure}{0.48\linewidth}
\begin{tikzpicture}[scale=0.6]
\pgfplotsset{
    tick label style={font=\small},
    label style={font=\small},
    legend style={font=\small},
    every axis/.append style={
    very thick,
    tick style={semithick}}
    }
\begin{axis}[
	title={\large $(n,m,r) = (25,50,10)$},
	xmin=-0.01,   xmax=1,
	ymin=-0.03,   ymax=1.03,
	xlabel={Relative remaining gap},
	ylabel={Fraction of instances},
	xtick={0,0.2,0.4,0.6,0.8,1},
	ytick={0,0.2,0.4,0.6,0.8,1}, 
	legend style={draw=none,
        at={(0.37,1)}},
    grid=major,
%   very thick,
    cycle list name=color list
        ]
\pgfplotstableread{size2.txt}\mydata;
\addplot table [x=poi,y=crmip] {\mydata};
\addplot table [x=poi,y=cr] {\mydata};
\addplot table [x=poi,y=mip] {\mydata};
\addplot table [x=poi,y=scip] {\mydata};
\addplot table [x=poi,y=mc] {\mydata};
 \legend{$ \texttt{CRMIP}$,$ \texttt{CR}$,$\texttt{MIP}/\texttt{MIP19}$, $\texttt{SCIP}$,$ \texttt{MC}$}
%\addplot[smooth, point meta=explicit] coordinates {
%	(-2,3) (-1.5,2) (-0.3,-0.2) 
%	(1,1.2) (2,2) (3,5)};
\end{axis}
\end{tikzpicture}
\end{subfigure}

\vspace{1em}

\begin{subfigure}{0.48\linewidth}
\begin{tikzpicture}[scale=0.6]
\pgfplotsset{
    tick label style={font=\small},
    label style={font=\small},
    legend style={font=\small},
    every axis/.append style={
    very thick,
    tick style={semithick}}
    }
\begin{axis}[
	title={\large $(n,m,r) = (50,100,20)$},
	xmin=-0.01,   xmax=1,
	ymin=-0.03,   ymax=1.03,
	xlabel={Relative remaining gap},
	ylabel={Fraction of instances},
	xtick={0,0.2,0.4,0.6,0.8,1},
	ytick={0,0.2,0.4,0.6,0.8,1}, 
	legend style={draw=none,
        at={(0.37,1)}},
    grid=major,
%   very thick,
    cycle list name=color list
        ]
\pgfplotstableread{size3.txt}\mydata;
\addplot table [x=poi,y=crmip] {\mydata};
\addplot table [x=poi,y=cr] {\mydata};
\addplot table [x=poi,y=mip] {\mydata};
\addplot table [x=poi,y=scip] {\mydata};
\addplot table [x=poi,y=mc] {\mydata};
 \legend{$ \texttt{CRMIP}$, $ \texttt{CR}$,$\texttt{MIP}/\texttt{MIP19}$, $\texttt{SCIP}$,$ \texttt{MC}$}
%\addplot[smooth, point meta=explicit] coordinates {
%	(-2,3) (-1.5,2) (-0.3,-0.2) 
%	(1,1.2) (2,2) (3,5)};
\end{axis}
\end{tikzpicture}
\end{subfigure}
\begin{subfigure}{0.48\linewidth}
\begin{tikzpicture}[scale=0.6]
\pgfplotsset{
    tick label style={font=\small},
    label style={font=\small},
    legend style={font=\small},
    every axis/.append style={
    very thick,
    tick style={semithick}}
    }
\begin{axis}[
	title={\large $(n,m,r) = (100,200,20)$},
	xmin=-0.01,   xmax=1,
	ymin=-0.03,   ymax=1.03,
	xlabel={Relative remaining gap},
	ylabel={Fraction of instances},
	xtick={0,0.2,0.4,0.6,0.8,1},
	ytick={0,0.2,0.4,0.6,0.8,1}, 
	legend style={draw=none,
        at={(0.37,1)}},
    grid=major,
%   very thick,
    cycle list name=color list
        ]
\pgfplotstableread{size4.txt}\mydata;
\addplot table [x=poi,y=crmip] {\mydata};
\addplot table [x=poi,y=cr] {\mydata};
\addplot table [x=poi,y=mip] {\mydata};
\addplot table [x=poi,y=scip] {\mydata};
\addplot table [x=poi,y=mc] {\mydata};
 \legend{$ \texttt{CRMIP}$, $ \texttt{CR}$,$\texttt{MIP}/\texttt{MIP19}$,$\texttt{SCIP}$,$ \texttt{MC}$}
%\addplot[smooth, point meta=explicit] coordinates {
%	(-2,3) (-1.5,2) (-0.3,-0.2) 
%	(1,1.2) (2,2) (3,5)};
\end{axis}
\end{tikzpicture}
\end{subfigure}
\caption{\revision{Relaxation cdfs, where the horizontal axis represents the relative remaining gap and the vertical axis represents the fraction of instances solved. The number of sub-gradient inequalities (resp. discretizations) is $5$ (resp. $1$).}}\label{fig:performance}
\end{figure}

We compare the strength of six relaxations.  For each instance $p$ in a set $\mathcal{P}$ of problems, we denote $v_{p,i}$ as the lower bound obtained by  solving the relaxation $i$, and denote $u_{p}$ as the upper bound using \textsc{Ipopt}~\cite{wachter2006implementation}. We report the estimate of percentage of gap that remains after solving each of the relaxations. This estimate uses the upper bound $u_p$ instead of the true optimal solution. In addition, we  use the weakest bound, $\min_i\{v_{p,i}\}$, as a reference for comparing the relative gaps. Then, we define $r_{p,i}$, the percentage of gap that remains after solving a relaxation, and $\mu_i(\alpha)$, for $0\le\alpha\le 1$ as
\[
r_{p,i} =  \frac{u_p - v_{p,i} }{u_p - \min_{i'}\{v_{p,i'}\}},\quad 
\mu_i(\alpha) = \frac{1}{\vert \mathcal{P} \vert} \bigl\vert \{ p \in \mathcal{P} \mid r_{p,i} \le \alpha \} \bigr\vert.
\]
Here, $\mu_i(\alpha)$ is the estimate of proportion of instances where relative gap is at least $\alpha$ for relaxation $i$. %The function $\mu_i(\alpha)$ is then the distribution function of the relative remaining gap for relaxation $i$. In other words, i
In Figure~\ref{fig:performance}, curves to the left correspond to stronger relaxations. The figure shows that \texttt{CRMIP} performs significantly better than any of the other relaxations. Surprisingly, \texttt{CR}, a linear programming relaxation, performs better than \texttt{MIP}. Table~\ref{table:computationtime} shows the average solution times of relaxations, \revision{where the solution time of \texttt{CR} and \texttt{CRMIP} includes the running time of the separation procedure.} \revision{Notice that \texttt{MIP} solves up to 2.5x faster than \texttt{MIP19} in our experiments.} These results demonstrate that relaxations in factorable-programming based solvers can be significantly improved by following the techniques proposed here. In Appendix~\ref{app:computation}, we report a computational study on how the quality of relaxations varies as the number of sub-gradient inequalities and discretization points change. 

\begin{table}[ht]
\begin{center}
\begin{tabular}{ | c |c | c |c | c | c| c |}
   \hline
 $(n,m,r)$  &  $\texttt{MC}$ & $\texttt{SCIP}$ & $\texttt{MIP19}$ & $\texttt{MIP}$ & $\texttt{CR}$  & $\texttt{CRMIP}$  \\ \hline
 $(15,30,10)$& 0.01 &0.52	& 0.17&	0.11&	1.54& 2.17		  \\ \hline
$(25,50,10)$&0.01& 0.83	&	0.34&0.25&	1.95		&	4.40   \\ \hline
$(50,100,20)$&0.02& 1.94&	1.70&0.97&1.61	&	6.87 	\\ \hline
$(100,200,20)$&0.04 & 3.74 &15.51 &	6.10&2.1&	29.86   \\ \hline
\end{tabular}      \caption{ \revision{Average solution times of relaxations on $100$ instances (seconds).}}~\label{table:computationtime}
\end{center}
\end{table}

\section{Extensions}\label{section:extension}
\subsection{Logarithmic model}
In Theorem~\ref{them:MIP-phi-z}, we show that~(\ref{eq:MIP-z}) yields an MILP formulation for~(\ref{eq:MIP-standard}). Interestingly, replacing the incremental model in~(\ref{eq:MIP-z}) by any formulation of $\{\Delta_H\}_{H\in \mathcal{H}}$ yields an MILP formulation for~(\ref{eq:MIP-standard}). Here, the ideality is established by invoking a decomposition result, which has been used in the literature~\cite{schrijver1983short,kim2022reciprocity}. 
\begin{lemma}\label{lemma:commonsimplex}
\revision{For $i  = 1,2$, let $D_i$ be a subset in the space of $(x_i,y)\in \R^{n_i+m}$,  and let $D :=\{(x,y) \mid (x_i,y) \in D_i \; i = 1,2 \} $. If $\proj_{y}(D_1) = \proj_{y}(D_2)$ is a finite set of affinely independent points then $\conv(D) = \{(x,y) \mid (x_i,y) \in \conv(D_i) \;  i =1,2\}$.}  
\end{lemma}
\begin{proposition}\label{prop:separateIntegral}
Let $S \subseteq \R^{d \times (n+1)} \times \Z^q $ be an MILP formulation of $\{\Delta_H\}_{H\in \mathcal{H}}$, where $q$ is a positive integer, and $D$ be its continuous relaxation. Then, an MIP formulation for~(\ref{eq:MIP-standard}) is given by
\begin{equation}\label{eq:separateIntegral}
\bigl\{(f,\phi,z,\gamma) \bigm|  \phi \leq \conc_\Delta (\phi \mcirc F)(z)  ,\  f = F(z),\ (z,\gamma) \in S \bigr\}.
\end{equation}
Moreover, this formulation is ideal if $\phi \mcirc F: \Delta \to \R $ is concave-extendable from $\vertex(\Delta)$  and \revision{$D$ can be expressed as $\bigl\{(z_1, \ldots, z_d, \gamma_1, \ldots, \gamma_d)\bigm| (z_i,\gamma_i) \in D_i \; \forall i \in [d]\bigr\}$, where $D_i$ is an integral polytope in $\R^{n+1} \times \R^{q_i}$ for each $i \in [d]$.}
\end{proposition}
\begin{compositeproof}
Let $E$ be the formulation in \eqref{eq:separateIntegral}, and $R$ be its continuous relaxation. The validity of $E$ follows from arguments similar to those in Theorem~\ref{them:MIP-phi-z}. \revision{Next, we prove the ideality of $E$. Let $C: = \{(f,\phi,z) \mid \phi \leq \conc_\Delta (\phi \mcirc F)(z)  ,\  f = F(z),\ z \in \Delta\}$. Let $R^0:=C$, and for $i \in [d]$, let $y^i :=(\gamma_1, \ldots, \gamma_i)$ and define $R^i := \{(f,\phi,z, y^i ) \mid (f,\phi,z,y^{i-1}) \in R^{i-1},\ (z_i,\gamma_i) \in D_i\}$ and 
\begin{equation*}
    B^i := \bigl\{(f,\phi,z,y^i) \in R^i \bigm| z \in \Z^{d \times(n+1)},\ \gamma_j \in \Z^{q_j} \; j =1, \ldots, i  \bigr\}.
\end{equation*} We will recursively invoke Lemma~\ref{lemma:commonsimplex} to argue that $R^d = \conv(B^d)$. This implies that $\vertex(R^d) \subseteq B^d$, proving the ideality of $E$ as $R^d = R$. By the concave extendability of $\phi \mcirc F$ from $\vertex(\Delta)$, we obtain that $R^0$ is the convex hull of $B^0:=\{(f,\phi,z) \in C \mid z \in \{0,1\}^{d \times(n+1)} \}$. Let $i \in [d]$, and suppose that $R^{i-1}= \conv(B^{i-1})$. Then,  
	\[
	\begin{aligned}
	R^i &= \Bigl\{(f,\phi,z,y^i) \Bigm| (f,\phi,z,y^{i-1}) \in \conv(B^{i-1}),\ (z_i,\gamma_i) \in \conv\bigl(D_i \cap (\Z^{n+1 +q_i})\bigr) \Bigr\} \\
	&= \conv\Bigl(\bigl\{(f,\phi,z,y^i) \bigm| (f,\phi,z,y^{i-1}) \in B^{i-1},\ (z_i,\gamma_i) \in D_i \cap \Z^{n+1+q_i} \bigr\} \Bigr),
	\end{aligned}
	\]
	where the first equality holds since $R^{i-1} = \conv(B^{i-1})$ and $D_i$ is integral, and the second equality follows from Lemma~\ref{lemma:commonsimplex} as $\proj_{z_i}(B^{i-1}) = \proj_{z_i}(D_i \cap \Z^{n+1+q_i} ) = \vertex(\Delta_i)$ is a set of affinely independent points. Therefore, we conclude that $R^d = \conv(B^d)$.}
	 \Halmos
\end{compositeproof}
\begin{remark}\label{rmk:log}
We can use the incremental formulation~\eqref{eq:Inc} as $S$ in Proposition~\ref{prop:separateIntegral} to obtain ideal formulations since \eqref{eq:Inc} is an integral formulation of $\{{\Delta_H}\}_{H\in \mathcal{H}}$. Instead, we can also choose $S$ so that it uses fewer binary variables. To simplify notation, we assume $l_i=n$ for $i \in [ d ]$. We select an interval using SOS2 constraints~\cite{beale1970special} on continuous variables $\lambda \in \Lambda_i:= \bigl\{\lambda \geq 0 \bigm| \sum_{j=0}^n \lambda_{ij}=1\bigr\}$, which allows at most two variables, $\lambda_{ij}$ and $\lambda_{ij'}$, to be non-zero and requires that $j$ and $j'$ are adjacent. For a sequence of distinct binary vectors $\{\eta^t\}_{t=1}^n \subseteq \{0,1\}^{\log_2(n)}$, where each adjacent pair $(\eta^{(t)}, \eta^{(t+1)})$ differ in at most one component,~\cite{vielma2011modeling} proposes the following ideal MILP formulation for $\lambda_i\in \Lambda_i$ satisfying SOS2 constraints, using binary variables, $\delta$, that are logarithmically many in the number of continuous variables $\lambda_i$:
\begin{equation}\label{eq:SOS2-log}
\begin{aligned}
&\lambda_i \in \Lambda_i \qquad \delta_{i} \in \{0,1\}^{\lceil \log_2(n) \rceil} \\
& \qquad \qquad \sum_{j \notin  L_{ik}} \lambda_{ij} \leq \delta_{ik} \leq 1-  \sum_{j \notin R_{ik}} \lambda_{ij}
 \text{ for } k  \in \bigl[ \lceil\log_2(n) \rceil \bigr],
\end{aligned}
\end{equation}
where $\eta^{(0)}:= \eta^{(1)}$ and $\eta^{(n+1)}: = \eta^{n}$, $L_{ik}: = \bigl\{j \in \{0\} \cup [  n ] \bigm| \eta^{(j)}_k = 1 \text{ or } \eta^{(j+1)}_k = 1  \bigr\}$, and $R_{ik}: = \bigl\{j \in \{0\} \cup [  n ] \bigm| \eta^{(j)}_k = 0 \text{ or } \eta^{(j+1)}_k = 0  \bigr\}$. Observe that integrality of \eqref{eq:SOS2-log} can also be inferred from Lemma~\ref{lemma:ideal} since $L_{ik}\cap R_{ik} = \emptyset$ and $\sum_{j\in L_{ik}\cup R_{ik}}\lambda_{ij}=1$ defines a face of $\Lambda_i$. Moreover, let $\Lambda = \prod_{i=1}^d\Lambda_i$ and define $T:\aff(\Delta)\rightarrow \aff(\Lambda)$ to be such that $T(z)_{ij} = z_{ij}-z_{ij+1}$, where $z_{in+1}$ is assumed to be zero. Then, $T$ is a transformation that maps $\Delta$ to $\Lambda$ and its inverse is given by $z_{ij} = \sum_{j'=j}^{n}\lambda_{ij'}$. We can rewrite \eqref{eq:SOS2-log} in $(z,\delta)$ space by simply replacing $\lambda_{ij}$ with $T(z)_{ij}$. It follows that the resulting set is integral. Therefore, if we use this set as $S$ in Proposition~\ref{prop:separateIntegral} we have that \eqref{eq:separateIntegral} is an ideal formulation whenever $(\phi\circ F)(\cdot)$ is concave extendable from $\vertex(\Delta)$. This representation requires $d\lceil\log_2(n)\rceil$ binary variables. \QED 
\end{remark}
\subsection{Multiple composite functions}
In this subsection, we generalize our results in Sections~\ref{section:MIP-conv} and~\ref{section:MIP-relaxation} to a vector of composite functions. Let $\kappa$ be a positive integer, and let $\theta \mcirc f: X \to \R^\kappa$ be a vector of composite functions defined as $(\theta \mcirc f)(x) = ((\theta_1 \mcirc f)(x), \ldots, (\theta_\kappa \mcirc f)(x))$, where the outer-function $\theta: \R^d \to \R^\kappa$ is defined as $\theta(f) = ( \theta_1(f), \ldots, \theta_\kappa(f))$ and the inner-function $f:X \to \R^d$ is defined as $f(x)=(f_1(x),\ldots,f_d(x))$.  
%Last, we generalize the constructions in Theorem~\ref{them:MIP-phi-z}  to the context of a vector of functions $\theta \mcirc f:X   \to \R^\kappa$ defined as $(\theta \mcirc f)(x) = \bigl((\theta_1 \mcirc f)(x), \ldots, (\theta_\kappa \mcirc f)(x)\bigr)$, where the outer-function $\theta: \R^d \to \R^\kappa$ is defined as $\theta(f) = \bigl( \theta_1(f), \ldots, \theta_\kappa(f)\bigr)$. More precisely, 
First, we generalize Theorem~\ref{them:MIP-phi-z} to model the following disjunctive constraints: 
\begin{equation}\label{eq:MIP-standard-simu}
	\bigcup_{H \in \mathcal{H}} \conv\bigl( \hypo(\theta|_H)\bigr),
\end{equation}
where $\hypo(\theta|_H): = \bigl\{(f,\theta) \mid \theta \leq \theta(f), f \in H\bigr\}$.  To obtain an MIP formulation of this constraint,  we convexify a vector of functions $\theta \mcirc F:\Delta \to \R^\kappa$ defined as $(\theta \mcirc F)(z) = \bigl( (\theta_1\mcirc F)(z), \ldots, (\theta_\kappa \mcirc F)(z)\bigr)$ over $\Delta$. 
\begin{proposition}\label{prop:MIP-phi-z-simu}
Let $S$ and $D$ be the sets defined as in Proposition~\ref{prop:separateIntegral}. An MIP formulation for~(\ref{eq:MIP-standard-simu}) is given by
\begin{equation}\label{eq:MIP-z-simu}
	\Bigl\{(f,\theta,z,\delta) \Bigm|  (z,\theta) \in  \conv\bigl( \hypo( \theta \mcirc F)\bigr)
 ,\ f = F(z),\ (z,\gamma) \in S \Bigr\},
\end{equation}
which is ideal if $\conv\bigl(\hypo(\theta \mcirc F|_{\vertex(\Delta)})\bigr) = \conv\bigl(\hypo(\theta \mcirc F)\bigr)$ and $D$ is integral. 
\end{proposition}
\begin{compositeproof}
The validity of~(\ref{eq:MIP-z-simu}) follows from the proof of Theorem~\ref{them:MIP-phi-z} by letting $E_C:= \bigl\{(f,\theta,z) \bigm| (z,\theta) \in \hypo(\theta\mcirc F),\ f = F(z),\ z \in C \bigr\}$ for $C \subseteq \Delta$. Its ideality  can be established using the same argument as in the proof of Proposition~\ref{prop:separateIntegral}. \Halmos 
\end{compositeproof}
 In general, individually convexifying the hypograph of $\theta_k \mcirc F$ does not yield the simultaneous convex hull of the hypograph of $\theta \mcirc F$. Thus, when variable  $z$ is restricted to a face $\Delta_H$ of $\Delta$, the formulation obtained by replacing $\conv( \hypo( \theta \mcirc F))$ with  $\cap_{k \in [\kappa]}\conv(\hypo(\theta_k \mcirc F))$ fails to model the disjunctive constraint~\eqref{eq:MIP-standard-simu}.  However,  when $\theta\mcirc F$ is a vector of concave-extendable supermodular functions it is shown in Corollary 7 of~\cite{he2022tractable} that $ \conv\bigl(\hypo(\theta\mcirc F) \bigr) = \cap_{k=1}^\kappa \conv\bigl(\hypo(\theta_k\mcirc F ) \bigr)$, where, for each $k \in [\kappa]$, the convex hull of $\hypo(\theta_k\mcirc F)$ is described in Proposition~\ref{prop:stair-Q}. This allows us to generalize Corollary~\ref{cor:MIP-phi-supermodular} to a vector of composite functions.

Next, we generalize Corollary~\ref{cor:mul} to the case where the outer-function is a vector of multilinear functions, \textit{i.e.}, for $k \in [ \kappa ] $, $\theta_k(s_{1n}, \ldots, s_{dn}) = \sum_{I \in \mathcal{I}_k} c^k_I\prod_{i \in I}s_{in} $, where $\mathcal{I}_k$ is a collection of subsets of $ [d]$.  To do so, we need to derive the convex hull of $\Theta^Q:=\{(s,\theta) \mid \theta = \theta(s_{1n}, \ldots,s_{dn} ),\ s\in Q \}$. Recall that Theorem~\ref{them:hull-Lambda} provides a convex hull description for the multilinear set $M^\Lambda$. Since $\Theta^Q$ can be obtained as an image of $M^\Lambda$ under an affine transformation, we obtain the convex hull of $\Theta^Q$ as follows:
\begin{equation}\label{eq:hull-Q} 
\begin{aligned}
 &\theta_k = \sum_{I \in \mathcal{I}_k}  c^k_I\sum_{p \in E}\biggl(\prod_{i \in I} a_{ip_i}\biggr) w_p \quad \text{ for } k  \in [ \kappa ],\\
 & w \geq 0 \quad \sum_{p \in E} w_p = 1 \quad  s_i = \sum_{j = 0}^{n} v_{ij} \sum_{p \in E:p_i = j}w_p \quad \text{  for } i \in [  d  ].
 \end{aligned}
\end{equation}

\begin{corollary}\label{cor:MIP-phi-mul}
Consider a vector of composite functions $\theta \mcirc f$, where $\theta: \R^d \to \R^\kappa$ be a vector of multilinear functions. An ideal formulation for $\{\conv(\graph(\theta|_H))\}_{H \in \mathcal{H}}$ is given by $\bigl\{(\theta, w,s,\delta) \bigm| (\theta, w, s) \in(\ref{eq:hull-Q}),\  Z(s) =z,\ (z, \delta) \in (\ref{eq:Inc-1}) \bigr\}$, which, together with  $u(x) \leq s$ and $(x,s_{\cdot n}) \in W$, yields an MICP relaxation for the graph of $\theta \mcirc f$.
\end{corollary}
\begin{compositeproof}
This result follows from Proposition~\ref{prop:MIP-phi-z-simu}, Remark~\ref{rmk:DCR-vector}, Corollary~\ref{eq:hull-Q} and the invertible transformation $Z$. \Halmos
\end{compositeproof}

%\input{computationplot.tex}

%
%\begin{figure}[ht]
%\includegraphics[scale=0.3]{size1.png}
%\includegraphics[scale=0.3]{size2.png}
%\includegraphics[scale=0.3]{size3.png}
%\includegraphics[scale=0.3]{size4.png}
%\end{figure}
%

\section{Conclusion}
In this paper, we developed new MIP relaxations for MINLPs. First, we derived ideal formulations using existing convexification results with ideal formulations of the faces of a specific simplotope. Second, we blended the above technique with recently developed composite relaxation technique by relating the variables used in composite relaxations with those describing the simplotope. Third, we obtained geometric insights into the source of improvement in relaxation quality and improved our relaxations by exploiting local bound information for underestimators. Fourth, we showed how factorable MINLPs can be converted to convex MINLPs using our techniques by propagating under- and over-estimators from expression trees. Finally, we showed that our relaxations significantly improve over the state-of-the-art relaxations. We close 60-70\% of the gap on several polynomial optimization instances relative to McCormick relaxations.

\appendix

\section{Appendix}~\label{app:proof}

\subsection{Proof of Corollary~\ref{cor:bilinear-stair}}\label{app:bilinear-stair}
Let $\bar{b}(s) = s_{1n}s_{2n}$ for $s \in Q$ and  let  $\min_{\pi \in \Pi}\check{b}^\pi(\cdot)$ and  $\min_{\pi \in \Pi}\hat{b}^\pi(\cdot)$ denote the function defined in~(\ref{eq:bilinear-stair-1}) and~(\ref{eq:bilinear-stair-2}), respectively. 
First, we use the constructions in Corollary 4 in~\cite{he2022tractable} to show that $\check{b}(\cdot)$  describes the convex envelope of $\bar{b}(\cdot)$ over $Q$. Let $\tilde{b}(s_1,s_2) := \bar{b}(s_1,\tilde{s}_2)$ for $s \in Q$, where $\tilde{s}_2$ is the image of an affine transformation defined as $\tilde{s}_{2j} = a_{20} + \sum_{k = 1}^j(a_{ik} - a_{ik-1})(1-z_{i,n+1-k})$ for $j \in \{0\} \cup [ n ]$, where $z$ denote $Z(s)$ defined as in~(\ref{eq:Z_trans}). Then, $\tilde{b}(\cdot)$ is concave extendable from $\vertex(Q)$, and, as $\tilde{b}(v_{1j_1},v_{1j_2}) = \bar{b}(v_{1j_1}, v_{2,n-j_2})$, it is submodular over $\vertex(Q)$. By Corollary 4 in~\cite{he2022tractable}, it suffices to verify that $\check{b}^\pi(\cdot)$ is obtained by affinely interpolating $\bar{b}(\cdot)$ over points $\bigl(v_{1p^0_1}, v_{2,n-p^0_2} \bigr), \bigl(v_{1p^1_1}, v_{2,n-p^1_2} \bigr), \ldots, \bigl(v_{1p^{2n}_1}, v_{2,n-p^{2n}_2}\bigr)$.
% \[    
% \bigl(v_{1p^0_1}, v_{2,n-p^0_2} \bigr), \bigl(v_{1p^1_1}, v_{2,n-p^1_2} \bigr), \ldots, \bigl(v_{1p^{2n}_1}, v_{2,n-p^{2n}_2}\bigr).
% \]
To see this, we observe that $\check{b}\bigl(v_{1p^0_1}, v_{2,n-p^0_2}\bigr) = a_{1 p^0_1} a_{2,n- p^0_2} = \check{b}\bigl(v_{1p^0_1}, v_{2p^0_2}\bigr)$, and for $t \in [ 2n ]$,
\[
\check{b}\bigl(v_{1p^t_1}, v_{2,n-p^t_2}  \bigr) = \check{b}\bigl(v_{1p^{t-1}_1}, v_{2,n-p^{t-1}_2}  \bigr) + \bar{b}\bigl(v_{1p^t_1}, v_{2,n-p^t_2}  \bigr) - \bar{b}\bigl(v_{1p^{t-1}_1}, v_{2,n-p^{t-1}_2}  \bigr). 
\]
In other words, $\check{b}\bigl(v_{1p^t_1}, v_{2,n-p^t_2}  \bigr) = \bar{b}\bigl(v_{1p^t_1}, v_{2,n-p^t_2}  \bigr)$ for $t \in [ 2n ]$. A similar argument can be used to show  that $\hat{b}(\cdot)$ describes the concave envelope of $\bar{b}(\cdot)$ over $Q$. More specifically, since the bilinear term $\bar{b}(\cdot)$ is concave-extendable from and supermodular over the vertices of $Q$,  by  Proposition~\ref{prop:stair-Q},  it suffices to verify that the affine function $\hat{b}^\pi(\cdot)$ is obtained by affinely interpolating $\bar{b}(\cdot)$ over $\vertex(\Upsilon_\pi)$. \QED 
%, it suffices to verify that the affine function $\hat{b}^\pi(\cdot)$ is obtained by affinely interpolating $\bar{b}(\cdot)$ over $\vertex(\Upsilon_\pi)$. To see this, we observe that $\hat{b}\bigl(v_{1p^0_1}, v_{2p^0_2}\bigr) = a_{1 p^0_1} a_{2 p^0_2} = \bar{b}\bigl(v_{1p^0_1}, v_{2p^0_2}\bigr)$, and for $t \in [ 2n ]$,
%\[
%\hat{b}\bigl(v_{1p^t_1}, v_{2p^t_2}  \bigr) = \hat{b}\bigl(v_{1p^{t-1}_1}, v_{2p^{t-1}_2}  \bigr) + \bar{b}\bigl(v_{1p^t_1}, v_{2p^t_2}  \bigr) - \bar{b}\bigl(v_{1p^{t-1}_1}, v_{2p^{t-1}_2}  \bigr). 
%\]
%In other words, $\hat{b}\bigl(v_{1p^t_1}, v_{2p^t_2}  \bigr) = \bar{b}\bigl(v_{1p^t_1}, v_{2p^t_2}  \bigr)$ for $t \in [ 2n ]$.
%
\subsection{Proof of Proposition~\ref{prop:eval}}\label{app:eval}
We first show that $\varphi_{\mathcal{H}}(\bar{x}, \bar{f}) = \conc_Q(\ephi)(s^*)$. By the first statement of Lemma~\ref{lemma:Inc} and the definition of $Z$, 
\[
\varphi_{\mathcal{H}}(\bar{x},\bar{f}) = \max\bigl\{ \conc_Q(\ephi)(s) \bigm| \u \leq s,\ s_{\cdot n} = \bar{f},\ s \in \cup_{H \in \mathcal{H}}Q_H \bigr\},
\]
where we denote the feasible region of the right-hand-side as $\mathcal{L}$. We will show that $s^*$ is the smallest point in $\mathcal{L}$, \textit{i.e.}, $s^* \in \mathcal{L}$ and, for every $s' \in \mathcal{L}$, $s^* \wedge s' = s^*$, where $\wedge$ denotes the component-wise minimum of two vectors. Then, $\max\{\conc_Q(\ephi)(s) \mid s \in \mathcal{L}\} \geq  \conc_Q(\ephi)(s^*) \geq \max\{\conc_Q(\ephi)(s) \mid s \in \mathcal{L}\}$, where the first inequality holds as $s^* \in \mathcal{L}$, and the second inequality holds because, by Lemma 8 in~\cite{he2021new}, $\conc_Q(\ephi)$ is non-increasing in $s_{ij}$ for all $i$ and $j \neq n$, and, for every point $s' \in \mathcal{L}$ we have $s' \geq s^*$ and $s'_{\cdot n} = s^*_{\cdot n}$. It follows that $\varphi_{\mathcal{H}}(\bar{x},\bar{f}) = \conc_Q(\phi)(s^*)$.

Now, we show that $s^*$ is the smallest point in $\mathcal{L}$. It follows from Proposition 4 in~\cite{he2021new} that $s^*_i \in  Q_i$ and $s^*_i \geq u^*_i$. As $u_i^* \geq \u_i$, we obtain that $s^*_i \geq \u_i$. Moreover, since $s^*_{ij}=a_{ij}$ for $j \leq \tau(i,\bar{t}_i-1)$ (resp. $s^*_{ij} = \bar{f}_i$ for $j > \tau(i,\bar{t}_i)$), it follows that $Z(s^*)_{ij} = 1$ for $j \leq \tau(i,\bar{t}_i-1)$ (resp. $Z(s^*)_{ij} = 0$ for $j > \tau(i,\bar{t}_i)$). In other words, $Z(s^*) \in \Delta_{\bar{H}}$ and, so, $s^* \in Q_{\bar{H}}$. Therefore, we conclude that $s^*\in \mathcal{L}$. Now, we show that, for $s' \in \mathcal{L}$, $s' \wedge s^* = s^*$. For any $s' \in \mathcal{L}$, we have $s'_{\cdot n}  \in \bar{H}$ and $Z(s') \in \cup_{H \in \mathcal{H}}\Delta_H$.
Then, by the second statement of Lemma~\ref{lemma:Inc} and $ s'_{\cdot n} = F\bigl(Z(s')\bigr)$, we obtain that $Z(s') \in \Delta_{\bar{H}}$.  Therefore, $s' \in Q_{\bar{H}}$, and in particular, $s'_{ij} = a_{ij}$ for $j \leq \tau(i,\bar{t}_i-1)$ and $s'_{ij} = \bar{f}_i$ for $j \geq \tau(i,\bar{t}_i)$.  This, together with $\u \leq s'$, implies that $ u^* \leq s'$, and thus $\xi_{i,a_i}(a;u^*_i) \leq \xi_{i,a_i}(a;s'_i) $ for every $a \in [a_{i0}, a_{in}]$. It turns out that for  $i \in \{1, \ldots, d\}$ and for  $j \in \{0, \ldots, n \}$, $(s'_i \wedge s^*_i)_j = \min\bigl\{ \conc(\xi_{i,a_i})(a_{ij};s'_i), \conc(\xi_{i,a_i})(a_{ij};u^*_i) \bigr\} = \conc(\xi_{i,a_i})(a_{ij};u^*_i) = s^*_{ij}$,
where the first equality holds since, by Lemma 5 in~\cite{he2021new}, for any $s_i \in Q_i$, $s_{ij} = \conc(\xi_{i,a_i})(a_{ij};s_i)$, and the second equality follows from $\xi_{i,a_i}(a;u^*_i) \leq \xi_{i,a_i}(a;s'_i)$. 

It is easy to see that $\varphi(\bar{x}, \bar{f}) = \conc_Q(\ephi) (\s)$ since $\bar{u}_i$ is the smallest vector larger than $u_i(\bar{x})$ with $\bar{u}_{in} = \bar{f}_i$. The proof that $\varphi_{\mathcal{H}-}(\bar{x}, \bar{f}) = \conc_Q(\ephi) (\hat{s})$ is similar to the proof above except that the feasible points are not required to satisfy $\bar{u}\le s$. Hence, the proof is complete since for $i\in [ d ]$, $s^*_i \geq \s_i$ and $s^*_i \geq \hat{s}_i$ with $s^*_{in} =\s_{in} = \hat{s}_{in}$, and by Lemma 8 in~\cite{he2021new},  $\conc_Q(\bar{\phi})$ is non-increasing in $s_{ij}$ for all $i$ and $j \neq n$. \QED
\subsection{Proof of Theorem~\ref{them:DCR+}}\label{app:DCR+}
Let $(x,f,\phi)$ be such that $\phi = \phi(f)$ and $f = f(x)$. We will construct $(s,z,\delta)$ such that $(x,\phi,s,z,\delta)$ satisfies the proposed relaxation. Let $t$ so that $f \in  \prod_{i=1}^d [a_{i \tau(i,t_i-1)},a_{i \tau(i,t_i)}]$. Let $\delta$ be binary so that for each $i$, $\delta_{ik} \geq \delta_{ik+1}$ and $1 = \delta_{it_i-1} > \delta_{it_i}=0$, and define $a'_{ij}=\BoundFunc{i:j:\delta_i}$ for all $j$. It follows from the  construction that for each $i$, $a'_{i0} = \cdots = a'_{i\tau(i,t_i-1)} \leq \cdots \leq a'_{i\tau(i,t_i)} = \cdots =a'_{in}$ and  there exists  $j'_i$ such that $s_{in} \in [a'_{ij'_i-1},a'_{ij'_i}]$. Then, we can choose $z_i$ to satisfy $1= z_{i0} = \cdots =z_{ij'_i-1} > z_{ij'_i}> z_{ij'_i+1} = \cdots =z_{in} =0$ and $f_i = G_{in}(z,\delta)$.  Now, for each $i$ and $j$, define $s_{ij} = G_{ij}(z,\delta)$, that is $s_{ij}=a'_{i0} + \sum_{j'=1}^j z_{ij'}(a'_{ij'}-a'_{ij'-1})$. It follows that $s_i = (a'_{i0}, \ldots, a'_{ij'_i-1}, f_i,\ldots, f_i)$. In other words, $s_{ij} = \max\{f_i,a'_{ij}\}$. Since $u_{ij}(x)$ underestimates $f(x)$ and is bounded from above by $a'_{ij}$, it follows that $u(x) \leq s$. Clearly, $(z,\delta)$ satisfies~(\ref{eq:Inc-1}), $(x,s_{\cdot}) \in W$, and $\phi \leq \phi(f) =  (\bar{\phi} \mcirc G)(z,\delta)$. Hence, the proof is complete.  \QED

\subsection{Random generation procedure}\label{app:generation}  
The instances are generated as follows:
\begin{itemize}
	\item each monomial $x^{\alpha_j}$ is generated by first randomly selecting nonzero entries of $\alpha_j$ such that  the number of nonzero entries is $2$ or $3$ with equal probability, then assign each nonzero entry $2$ or $3$ uniformly. 
	\item each entry of $d$ and $B$ is zero with probability $0.3$ and uniformly generated from $[0,1]$ with probability $0.7$.
	\item each entry of $A$ is uniformly distributed over $[-10,10]$, and $x^L_i$ and $x_i^U$ is uniformly selected from $\{0,1,2\}$ and $\{3,4\}$, respectively.
	\item $c = \sum_{j = 1}^m \nabla m_j(\tilde{x}) d_j$, and $b = A\tilde{x} + B\tilde{y}$, where $\tilde{x}$ is randomly generated from $[x^L, x^U]$, $\nabla m_j(\tilde{x})$ is the gradient of $x^{\alpha_j}$ at $\tilde{x}$, and $\tilde{y} = (\tilde{x}^{\alpha_1}, \ldots, \tilde{x}^{\alpha_m} )$.
\end{itemize}

\subsection{Additional computational results on the quality of relaxations}~\label{app:computation}
\begin{figure}[h]
\begin{subfigure}{0.48\linewidth}
\begin{tikzpicture}[scale=0.6]
\pgfplotsset{
    tick label style={font=\small},
    label style={font=\small},
    legend style={font=\small},
    every axis/.append style={
    very thick,
    tick style={semithick}}
    }
\begin{axis}[
	xmin=-0.01,   xmax=1,
	ymin=-0.03,   ymax=1.03,
	xlabel={Relative remaining gap},
	ylabel={Fraction of instances},
	xtick={0,0.2,0.4,0.6,0.8,1},
	ytick={0,0.2,0.4,0.6,0.8,1}, 
	legend style={draw=none,
        at={(0.97,0.7)}},
    grid=major,
%   very thick,
    cycle list name=color list
        ]
\pgfplotstableread{size3plus.txt}\mydata;
\addplot table [x=poi,y=mc] {\mydata};
\addplot table [x=poi,y=cr] {\mydata};
\addplot table [x=poi,y=cr1] {\mydata};
\addplot table [x=poi,y=cr3] {\mydata};
%\addplot table [x=poi,y=cr2] {\mydata};
 \legend{$ \texttt{MC}$, $ \texttt{CR(5)}$,$\texttt{CR(7)}$,$\texttt{CR(9)}$}
%\addplot[smooth, point meta=explicit] coordinates {
%	(-2,3) (-1.5,2) (-0.3,-0.2) 
%	(1,1.2) (2,2) (3,5)};
\end{axis}
\end{tikzpicture}
\end{subfigure}
\begin{subfigure}{0.48\linewidth}
\begin{tikzpicture}[scale=0.6]
\pgfplotsset{
    tick label style={font=\small},
    label style={font=\small},
    legend style={font=\small},
    every axis/.append style={
    very thick,
    tick style={semithick}}
    }
\begin{axis}[
	xmin=-0.01,   xmax=1,
	ymin=-0.03,   ymax=1.03,
	xlabel={Relative remaining gap},
	ylabel={Fraction of instances},
	xtick={0,0.2,0.4,0.6,0.8,1},
	ytick={0,0.2,0.4,0.6,0.8,1}, 
	legend style={draw=none,
        at={(0.97,0.7)}},
    grid=major,
%   very thick,
    cycle list name=color list
        ]
\pgfplotstableread{size3plus.txt}\mydata;
\addplot table [x=poi,y=mc] {\mydata};
\addplot table [x=poi,y=crmip] {\mydata};
\addplot table [x=poi,y=crmip1] {\mydata};
\addplot table [x=poi,y=crmip2] {\mydata};
\addplot table [x=poi,y=crmip3] {\mydata};
 \legend{$ \texttt{MC}$, $ \texttt{CRMIP(5,1)}$,$\texttt{CRMIP(7,1)}$,$\texttt{CRMIP(7,2)}$,$\texttt{CRMIP(9,2)}$}
%\addplot[smooth, point meta=explicit] coordinates {
%	(-2,3) (-1.5,2) (-0.3,-0.2) 
%	(1,1.2) (2,2) (3,5)};
\end{axis}
\end{tikzpicture}
\end{subfigure}
\caption{Relaxation cdfs comparing relaxation strengths on $100$ instances of size $(n,m,r) = (50,100,20)$.}\label{fig:performance3}
\end{figure}
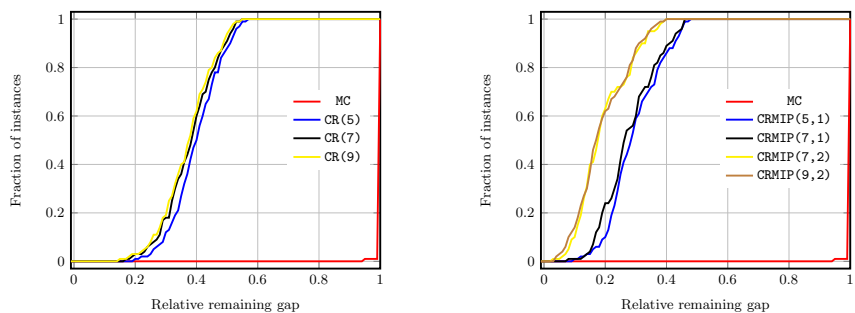

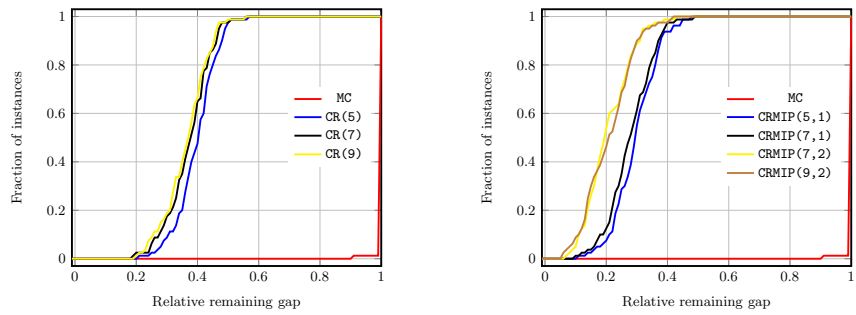
\begin{figure}[h]
\begin{subfigure}{0.48\linewidth}
\begin{tikzpicture}[scale=0.6]
\pgfplotsset{
    tick label style={font=\small},
    label style={font=\small},
    legend style={font=\small},
    every axis/.append style={
    very thick,
    tick style={semithick}}
    }
\begin{axis}[
	xmin=-0.01,   xmax=1,
	ymin=-0.03,   ymax=1.03,
	xlabel={Relative remaining gap},
	ylabel={Fraction of instances},
	xtick={0,0.2,0.4,0.6,0.8,1},
	ytick={0,0.2,0.4,0.6,0.8,1}, 
	legend style={draw=none,
        at={(0.97,0.7)}},
    grid=major,
%   very thick,
    cycle list name=color list
        ]
\pgfplotstableread{size4plus.txt}\mydata;
\addplot table [x=poi,y=mc] {\mydata};
\addplot table [x=poi,y=cr] {\mydata};
\addplot table [x=poi,y=cr1] {\mydata};
\addplot table [x=poi,y=cr3] {\mydata};
%\addplot table [x=poi,y=cr2] {\mydata};
 \legend{$ \texttt{MC}$, $ \texttt{CR(5)}$,$\texttt{CR(7)}$,$\texttt{CR(9)}$}
%\addplot[smooth, point meta=explicit] coordinates {
%	(-2,3) (-1.5,2) (-0.3,-0.2) 
%	(1,1.2) (2,2) (3,5)};
\end{axis}
\end{tikzpicture}
\end{subfigure}
\begin{subfigure}{0.48\linewidth}
\begin{tikzpicture}[scale=0.6]
\pgfplotsset{
    tick label style={font=\small},
    label style={font=\small},
    legend style={font=\small},
    every axis/.append style={
    very thick,
    tick style={semithick}}
    }
\begin{axis}[
	xmin=-0.01,   xmax=1,
	ymin=-0.03,   ymax=1.03,
	xlabel={Relative remaining gap},
	ylabel={Fraction of instances},
	xtick={0,0.2,0.4,0.6,0.8,1},
	ytick={0,0.2,0.4,0.6,0.8,1}, 
	legend style={draw=none,
        at={(0.97,0.7)}},
    grid=major,
%   very thick,
    cycle list name=color list
        ]
\pgfplotstableread{size4plus.txt}\mydata;
\addplot table [x=poi,y=mc] {\mydata};
\addplot table [x=poi,y=crmip] {\mydata};
\addplot table [x=poi,y=crmip1] {\mydata};
\addplot table [x=poi,y=crmip2] {\mydata};
\addplot table [x=poi,y=crmip3] {\mydata};
 \legend{$ \texttt{MC}$, $ \texttt{CRMIP(5,1)}$,$\texttt{CRMIP(7,1)}$,$\texttt{CRMIP(7,2)}$,$\texttt{CRMIP(9,2)}$}
%\addplot[smooth, point meta=explicit] coordinates {
%	(-2,3) (-1.5,2) (-0.3,-0.2) 
%	(1,1.2) (2,2) (3,5)};
\end{axis}
\end{tikzpicture}
\end{subfigure}
\caption{Relaxation cdfs comparing relaxation strengths on $100$ instances of size $(n,m,r) = (100,200,20)$.}\label{fig:performance4}
\end{figure}

Here, we investigate how the quality of relaxations varies as the number of sub-gradient inequalities, denoted as $\alpha$, and the number of discretization points, denoted as $\beta$, increases. 
\begin{itemize}
	\item \texttt{CR($\alpha$)}: This relaxation is obtained by increasing the number of sub-gradient inequalities used in the $\texttt{CR}$ from $5$ to $\alpha$.
	\item \texttt{CRMIP($\alpha,\beta$)}: This relaxation is obtained by  changing the number of sub-gradient inequalities and of discretization points used in $\texttt{CRMIP}$ from $(5,1)$ to $(\alpha,\beta)$. Here, we limit the solution time to 3600 seconds.
\end{itemize}
Figures~\ref{fig:performance1},~\ref{fig:performance2},~\ref{fig:performance3} and~\ref{fig:performance4} plot relaxation cdfs comparing \texttt{MC} with relaxations \texttt{CR($\alpha$)}  and \texttt{CRMIP($\alpha,\beta$)} on various settings of parameters $(n,m,r)$. Tables~\ref{table:computationtime2} and~\ref{table:computationtime3} shows the average solution times of these settings.

\begin{figure}[hp]
\begin{subfigure}{0.48\linewidth}
\begin{tikzpicture}[scale=0.6]
\pgfplotsset{
    tick label style={font=\small},
    label style={font=\small},
    legend style={font=\small},
    every axis/.append style={
    very thick,
    tick style={semithick}}
    }
\begin{axis}[
	xmin=-0.01,   xmax=1,
	ymin=-0.03,   ymax=1.03,
	xlabel={Relative remaining gap},
	ylabel={Fraction of instances},
	xtick={0,0.2,0.4,0.6,0.8,1},
	ytick={0,0.2,0.4,0.6,0.8,1}, 
	legend style={draw=none,
        at={(0.97,0.7)}},
    grid=major,
%   very thick,
    cycle list name=color list
        ]
\pgfplotstableread{size1plus.txt}\mydata;
\addplot table [x=poi,y=mc] {\mydata};
\addplot table [x=poi,y=cr] {\mydata};
\addplot table [x=poi,y=cr1] {\mydata};
\addplot table [x=poi,y=cr3] {\mydata};
%\addplot table [x=poi,y=cr2] {\mydata};
 \legend{$ \texttt{MC}$, $ \texttt{CR(5)}$,$\texttt{CR(7)}$,$\texttt{CR(9)}$}
%\addplot[smooth, point meta=explicit] coordinates {
%	(-2,3) (-1.5,2) (-0.3,-0.2) 
%	(1,1.2) (2,2) (3,5)};
\end{axis}
\end{tikzpicture}
\end{subfigure}
\begin{subfigure}{0.48\linewidth}
\begin{tikzpicture}[scale=0.6]
\pgfplotsset{
    tick label style={font=\small},
    label style={font=\small},
    legend style={font=\small},
    every axis/.append style={
    very thick,
    tick style={semithick}}
    }
\begin{axis}[
	xmin=-0.01,   xmax=1,
	ymin=-0.03,   ymax=1.03,
	xlabel={Relative remaining gap},
	ylabel={Fraction of instances},
	xtick={0,0.2,0.4,0.6,0.8,1},
	ytick={0,0.2,0.4,0.6,0.8,1}, 
	legend style={draw=none,
        at={(0.97,0.7)}},
    grid=major,
%   very thick,
    cycle list name=color list
        ]
\pgfplotstableread{size1plus.txt}\mydata;
\addplot table [x=poi,y=mc] {\mydata};
\addplot table [x=poi,y=crmip] {\mydata};
\addplot table [x=poi,y=crmip1] {\mydata};
\addplot table [x=poi,y=crmip2] {\mydata};
\addplot table [x=poi,y=crmip3] {\mydata};
 \legend{$ \texttt{MC}$, $ \texttt{CRMIP(5,1)}$,$\texttt{CRMIP(7,1)}$,$\texttt{CRMIP(7,2)}$,$\texttt{CRMIP(9,2)}$}
%\addplot[smooth, point meta=explicit] coordinates {
%	(-2,3) (-1.5,2) (-0.3,-0.2) 
%	(1,1.2) (2,2) (3,5)};
\end{axis}
\end{tikzpicture}
\end{subfigure}
\caption{Relaxation cdfs comparing relaxation strengths on $100$ instances of size $(n,m,r) = (15,30,10)$.}\label{fig:performance1}
\end{figure}
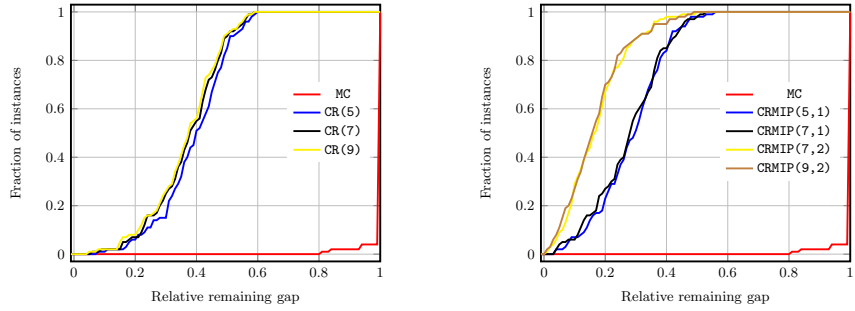

\begin{figure}[h]
\begin{subfigure}{0.48\linewidth}
\begin{tikzpicture}[scale=0.6]
\pgfplotsset{
    tick label style={font=\small},
    label style={font=\small},
    legend style={font=\small},
    every axis/.append style={
    very thick,
    tick style={semithick}}
    }
\begin{axis}[
	xmin=-0.01,   xmax=1,
	ymin=-0.03,   ymax=1.03,
	xlabel={Relative remaining gap},
	ylabel={Fraction of instances},
	xtick={0,0.2,0.4,0.6,0.8,1},
	ytick={0,0.2,0.4,0.6,0.8,1}, 
	legend style={draw=none,
        at={(0.97,0.7)}},
    grid=major,
%   very thick,
    cycle list name=color list
        ]
\pgfplotstableread{size2plus.txt}\mydata;
\addplot table [x=poi,y=mc] {\mydata};
\addplot table [x=poi,y=cr] {\mydata};
\addplot table [x=poi,y=cr1] {\mydata};
\addplot table [x=poi,y=cr3] {\mydata};
%\addplot table [x=poi,y=cr2] {\mydata};
 \legend{$ \texttt{MC}$, $ \texttt{CR(5)}$,$\texttt{CR(7)}$,$\texttt{CR(9)}$}
%\addplot[smooth, point meta=explicit] coordinates {
%	(-2,3) (-1.5,2) (-0.3,-0.2) 
%	(1,1.2) (2,2) (3,5)};
\end{axis}
\end{tikzpicture}
\end{subfigure}
\begin{subfigure}{0.48\linewidth}
\begin{tikzpicture}[scale=0.6]
\pgfplotsset{
    tick label style={font=\small},
    label style={font=\small},
    legend style={font=\small},
    every axis/.append style={
    very thick,
    tick style={semithick}}
    }
\begin{axis}[
	xmin=-0.01,   xmax=1,
	ymin=-0.03,   ymax=1.03,
	xlabel={Relative remaining gap},
	ylabel={Fraction of instances},
	xtick={0,0.2,0.4,0.6,0.8,1},
	ytick={0,0.2,0.4,0.6,0.8,1}, 
	legend style={draw=none,
        at={(0.97,0.7)}},
    grid=major,
%   very thick,
    cycle list name=color list
        ]
\pgfplotstableread{size2plus.txt}\mydata;
\addplot table [x=poi,y=mc] {\mydata};
\addplot table [x=poi,y=crmip] {\mydata};
\addplot table [x=poi,y=crmip1] {\mydata};
\addplot table [x=poi,y=crmip2] {\mydata};
\addplot table [x=poi,y=crmip3] {\mydata};
 \legend{$ \texttt{MC}$, $ \texttt{CRMIP(5,1)}$,$\texttt{CRMIP(7,1)}$,$\texttt{CRMIP(7,2)}$,$\texttt{CRMIP(9,2)}$}
%\addplot[smooth, point meta=explicit] coordinates {
%	(-2,3) (-1.5,2) (-0.3,-0.2) 
%	(1,1.2) (2,2) (3,5)};
\end{axis}
\end{tikzpicture}
\end{subfigure}
\caption{Relaxation cdfs comparing relaxation strengths on $100$ instances of size $(n,m,r) = (25,50,10)$.}\label{fig:performance2}
\end{figure}
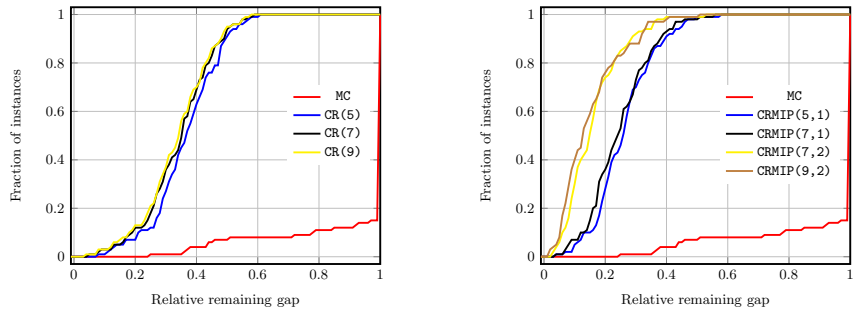

\begin{table}[ht]
\begin{center}
\begin{tabular}{ | c |c | c | c| c |}
    \hline
 $(n,m,r)$  & $\texttt{CR(5)}$ & $\texttt{CR(7)}$ & $\texttt{CR(9)}$   \\ \hline
 $(15,30,10)$ &1.54	&	1.49&	1.79  \\ \hline
$(25,50,10)$& 1.95	&	1.54&	2.53		 \\ \hline
$(50,100,20)$ & 1.61 &	2.32&5.32		\\ \hline
$(100,200,20)$ & 2.1 &	8.03&18.32		\\ \hline
\end{tabular}      \caption{Average solution times of relaxations \texttt{CR($\alpha$)} on $100$ instances (seconds).}~\label{table:computationtime2}
\end{center}
\end{table}

\begin{table}[ht]
\begin{center}
\begin{tabular}{ | c  |c | c | c| c |}
    \hline
 $(n,m,r)$   & $\texttt{CRMIP(5,1)}$ & $\texttt{CRMIP(7,1)}$ & $\texttt{CRMIP(7,2)}$  & $\texttt{CRMIP(9,2)}$  \\ \hline
 $(15,30,10)$  &1.91	&	3.5&	4.72& 8.56		  \\ \hline
$(25,50,10)$ & 4.12	&	6.27&	9.24		&	21.12   \\ \hline
$(50,100,20)$ & 6.87 &	20.5&57.95	&	169.31 	\\ \hline
$(100,200,20)$ & 29.86 &	153&1065.76	&	2367 	\\ \hline
\end{tabular}      \caption{Average solution times of relaxations \texttt{CRMIP($\alpha, \beta$)} on $100$ instances (seconds).}\label{table:computationtime3}
\end{center}
\end{table}

%\section*{Acknowledgments}
%We would like to acknowledge the assistance of volunteers in putting
%together this example manuscript and supplement.

\bibliographystyle{siamplain}
\bibliography{reference}

\end{document}